\documentclass[11pt]{article}

\usepackage{graphicx}
\usepackage{amsmath}
\usepackage{amsfonts}
\usepackage{fancyhdr} 
\usepackage{mathrsfs}

\usepackage{wrapfig}
\usepackage{bbold}

\usepackage{float}

\newtheorem{thm}{Theorem}
\newtheorem{lem}{Lemma}

\newtheorem{prop}[lem]{Proposition}

\numberwithin{equation}{section}
\numberwithin{defn}{section}
\newcommand{\real}{\mathbb{R}}

\newcommand{\integer}{\mathbb{Z}}

\newcommand{\imp}{\mathscr{I}}
\newcommand{\pimp}{\mathscr{P}}

\newenvironment{pf}{\noindent {\em Proof}.\ \ }{\hspace*{\fill}\rule{.5ex}{1.4ex}\,}

\title{\vspace*{-1in}The refined impedance transform for 1D acoustic reflection data}%Acoustic impedance across a layered slab
\author{Peter C.~Gibson\footnote{Dept.~of Mathematics \& Statistics, York University, 4700 Keele St., Toronto, Ontario, Canada, M3J~1P3, $\mathtt{pcgibson@yorku.ca}$} }
\date{January 31, 2018}
\makeatletter
\let\newtitle\@title
\let\newauthor\@author
\let\newdate\@date
\makeatother

\begin{document}

\maketitle

\begin{abstract}

The one dimensional wave equation serves as a basic model for imaging modalities such as seismic which utilize acoustic data reflected back from a layered medium.  In 1955 Peterson et al.~described a  single scattering approximation for the one dimensional wave equation that relates the reflection Green's function to acoustic impedance.  The approximation is simple, fast to compute and has become a standard part of seismic theory.  The present paper re-examines this classical approximation in light of new results concerning the (exact) measurement operator for reflection imaging of layered media, and shows that the classical approximation can be substantially improved.  We derive an alternate formula, called the refined impedance transform, that retains the simplicity and speed of computation of the classical estimate, but which is qualitatively more accurate and applicable to a wider range of recorded data.  The refined impedance transform can be applied to recorded data directly (without the need to deconvolve the source wavelet), and solves exactly the inverse problem of determining the value of acoustic impedance on the far side of an arbitrary slab of unknown structure.  The results are illustrated with numerical examples.  
%Numerical examples illustrate the results.  
\end{abstract}
\begin{center}MSC 35L05, 35Q86, 35R30; Keywords: one dimensional wave equation, impedance inversion, inverse problems\end{center}

%MSC 35L05, 35R30

\section{Introduction}
We consider the one dimensional wave equation governing particle velocity $u(x,t)$:
%in a layered acoustic medium
\begin{subequations}\label{wave}
\begin{gather}
u_{tt}-\textstyle\frac{1}{\zeta}(\zeta u_x)_x=0\label{wave-equation}\\
u(x,0)=W(x)\qquad u_t(x,0)=-W^\prime(x)\label{initial-condition}
\end{gather}
\end{subequations}
where $\zeta(x)$ denotes acoustic impedance as a function of position (in units of one way travel time) and $W$ denotes a compactly supported initial waveform.  The impedance function $\zeta$ is assumed to be positive, bounded above and bounded away from 0, and to have respective constant values $\zeta_-$ and $\zeta_+$ on intervals $(-\infty,x_-)$ and $[x_+,\infty)$, for some $0<x_-<x_+$.  Furthermore, the compact support of $W$ is assumed to be contained in the interval $(-\infty,x_-)$.  The reflection Green's function for (\ref{wave-equation}) is defined to be the restriction to $x=0$ of the unique solution $u(x,t)$ to (\ref{wave}) in the case where $W=\delta$.  More precisely,
\begin{equation}\label{Green}
G_\zeta(t)=\left\{\begin{array}{cc}
0&\mbox{ if }t\leq 0\\
u(0,t)&\mbox{ if }t>0
\end{array}
\right..
\end{equation}
The reflection data for (\ref{wave}), denoted $D(t)$, is defined to be the left-moving part of the unique solution for general $W$ at the spatial location $x=0$.  Thus
\begin{equation}\label{data}
D=\widetilde{W}\ast G_\zeta
\end{equation}
where $\widetilde{W}(x)=W(-x)$ is the reversal of $W$.   

The equation for pressure, as opposed to particle velocity, is obtained from (\ref{wave-equation}) by replacing $\zeta$ with its reciprocal.  The results for pressure are essentially the same as those for velocity, with the typical difference being a change in sign, since $G_{1/\zeta}=-G_\zeta$.   Thus replacing $G_\zeta$ with its negative in the various formulas below converts from velocity data to pressure data; see Appendix~\ref{sec-pressure}. 

\subsection{Background}

For consistency of notation we briefly derive the classical estimate of Peterson et al.~\cite{PeFiCo:1955} in terms of equation (\ref{wave}) so as to obtain the formulas appropriate to particle velocity $u(x,t)$ expressed in terms of one way travel time $x$ (see also \cite{OlScLe:1983}).   
Let $P=\{x_1,\ldots,x_{n+1}\}$ be a partition of $[a,b]$, so that 
$
a=x_1<\cdots<x_{n+1}=b
$,
and set $\Delta=\max_{1\leq j\leq n}x_{j+1}-x_j$.  Set $x^\ast_j=(x_j+x_{j+1})/2$ for $1\leq j\leq n$, and define the step function approximation $\zeta^P$ to $\zeta$ by the formula
\begin{equation}\label{step}
\zeta^P(x)=\left\{\begin{array}{cc}
\zeta(x_1)&\mbox{ if }x<x_1^\ast\\
\zeta(x_j)&\mbox{ if }x_{j-1}^\ast\leq x<x_j^\ast\quad (2\leq j\leq n)\\
\zeta(x_{n+1})&\mbox{ if }x_n^\ast\leq x
\end{array}\right..
\end{equation}
Thus $\zeta^P$ has $n$ jump points $x_j^\ast$ $(1\leq j\leq n)$ and agrees with $\zeta$ on $P$.  The reflectivities corresponding to these jump points are 
\begin{equation}\label{reflectivities}
r_j=\frac{\zeta(x_j)-\zeta(x_{j+1})}{\zeta(x_j)+\zeta(x_{j+1})}\qquad(1\leq j\leq n),
\end{equation}
and the reflectivity distribution corresponding to $\zeta^P$, expressed in terms of 2-way travel time, is 
\begin{equation}\label{reflectivity-distribution}
r^P(t)=\sum_{j=1}^nr_j\delta(t-2x_j^\ast).
\end{equation}
In the case where the partition $P$ is evenly spaced, the initial integral of $r^P$ up to a given time $\tau\in[2x_j^\ast,2x_{j+1}^\ast)$ for $1\leq j\leq n-1$, or $\tau\geq 2x_j^\ast=2x_n^\ast$, is 
\begin{equation}\label{partial-integral}
\int_{-\infty}^{\tau}r^P(t)\,dt=\int_{-\infty}^{2x_{j+1}}\frac{1}{2\Delta}\sum_{k=1}^{n}r_k\chi_{[2x_k,2x_{k+1})}(t)\,dt.
\end{equation}
The reflectivity function $R_\zeta$ determined by $\zeta$ is defined to be the pointwise limit of the integrand on the right-hand side of (\ref{partial-integral}) with fixed endpoints $a\leq x_-$ and $b\geq x_+$, 
\begin{equation}\label{limit}
\lim_{\Delta\rightarrow0}\frac{1}{2\Delta}\sum_{k=1}^{n}r_k\chi_{[2x_k,2x_{k+1})}(t)=\lim_{\Delta\rightarrow0}\sum_{k=1}^{n}\frac{\zeta(x_j)-\zeta(x_{j+1})}{2\Delta(\zeta(x_j)+\zeta(x_{j+1}))}\chi_{[2x_k,2x_{k+1})}(t)=-\frac{\zeta^\prime(t/2)}{4\zeta(t/2)}. 
\end{equation}
Thus
\begin{equation}\label{R}
R_\zeta(t):=-\frac{\zeta^\prime(t/2)}{4\zeta(t/2)}\qquad(t\in\real).
\end{equation}
It follows that, provided $\zeta$ is absolutely continuous, it may be expressed in terms of the reflectivity function by the equation
\begin{equation}\label{standard-inverse}
\zeta(x)=\zeta_-e^{-2\int_{-\infty}^{2x}R_\zeta(t)\,dt}.
\end{equation}
The standard formula attributed to Peterson et al.~rests on the approximation
\begin{equation}\label{bad-approximation}
G_\zeta\cong R_\zeta 
\end{equation}
which, using (\ref{standard-inverse}), yields the estimate
\begin{equation}\label{standard}
\zeta(x)\cong\zeta_-e^{-2\int_{-\infty}^{2x}G_\zeta(t)\,dt}.
\end{equation}
\pagestyle{fancyplain}
The approximation (\ref{bad-approximation}) only takes account of single scattering, and is therefore not exact.  Yet versions of the formula (\ref{standard}) continue to serve as a point of reference in the geophysics literature more than sixty years after its initial publication (see, e.g., \cite{TaMc:2017}, \cite{GhMa:2016} and references therein). There are several reasons for this: (i) the formula (\ref{standard}) is simple, and fast to evaluate given (properly deconvolved) measured data; (ii) the formula, while inexact, seems reasonably accurate for a certain limited class of impedance functions; and (iii) in the absence of explicit formulas for $G_\zeta$, it's not clear how (\ref{bad-approximation}) should be modified.   

The present paper introduces a new alternative to (\ref{standard}) which is both simpler and more accurate---optimal among formulas of the form $f\circ\int_{-\infty}^{2x}G_\zeta$ for some function $f$.  Recent insight into measurement operator $\zeta\mapsto D=\widetilde{W}\ast G_\zeta$ has made this improvement possible by eliminating reliance on the single scattering approximation (\ref{bad-approximation}).  The new refined impedance transform (defined below in \S\ref{sec-main}) allows one to solve exactly the inverse problem of determining acoustic impedance on the far side of an arbitrary layered slab, and gives a remarkably accurate approximation to the full impedance profile in a wider range of cases than the standard approximation (\ref{standard}).  

\subsection{Related literature\label{sec-related}}

Variants of (\ref{wave}) have been studied extensively over the last half century because they serve as a basic model in acoustic and electromagnetic imaging.  Much of what is known places certain restrictions on $\zeta$, such as that its logarithmic derivative $\zeta^\prime/\zeta$ be in $L^2(\real)$ \cite{SaSy:1988},\cite{SyWi:1998}, which implies that $\zeta$ is at least continuous.  Continuity of $\zeta$ (or smoothness) is an assumption behind much of the classical theory in the subject \cite{Ne:1981}.  An alternative approach is to discretize $\zeta$ on an even grid according to the Goupillaud scheme \cite{Go:1961},\cite{BuBu:1983}. 

Whatever the restrictions on $\zeta$, the various schemes for its determination from reflection data, including the textbook seismic approach \cite{BlCoSt:2001}, \cite{Yi:2001}, involve either approximations or iterative procedures or both.  Apart from (\ref{standard}), there is no straightforward transform analogous to the Radon transform of CT that directly converts the data into something representing physical structure.  

Recent results concerning a class of hyperbolic equations that includes (\ref{wave}) allow $\zeta$ to be completely general (see \cite{BlStSy:2013} and \cite{KiRi:2014}), but these are focussed on Newton-like schemes that require an initial guess and subsequent iteration.  Random media \cite{FoGaPaSo:2007} is another perspective from which (\ref{wave}) has been thoroughly explored albeit nondeterministically.  None of this gives a simple direct method to recover the coefficient $\zeta$ in (\ref{wave}) from time limited reflection data.  

The present paper details several features of the refined impedance transform that collectively set it apart from established results:
\begin{itemize}
\item it is simple, and fast to compute---fast enough to convert recorded data in real time;
\item the coefficient $\zeta$ may be essentially arbitrary (see \S\ref{sec-derivation} below);  
\item it recovers impedance on the far side of an unknown slab exactly, without deconvolution.
\end{itemize}

%\subsection{New results on scattering}
\subsection{Main result\label{sec-main}}
Define the \emph{refined impedance transform} of an integrable distribution $g\in\mathcal{S}^\prime(\real)$, depending on a pair of real parameters $w\neq0$ and $c>0$, by the formula
\begin{equation}\label{impedance-transform}
\mathscr{I}_{w,c}\;g(x)=c\frac{w-\int_{-\infty}^{2x}g(t)\,dt}{w+\int_{-\infty}^{2x}g(t)\,dt}.
\end{equation}
\begin{thm}\label{thm-slab}
Let $D$ denote the reflection data for (\ref{wave}) as defined by (\ref{data}), and set $w=\int W$.  If $w\neq 0$, then 
\begin{equation}\label{convergence}
\imp_{w,\zeta_-}D(x)\rightarrow\zeta_+\qquad\mbox{ as }\qquad x\rightarrow\infty.
\end{equation}
\end{thm}
This means that for any source wavelet which is not zero mean, the impedance value $\zeta_+$ on the far side of an arbitrary slab $[x_-,x_+]$ can in principle be recovered with arbitrary accuracy from reflection data recorded at $x=0$, provided the recording duration $T>2x_+$ is sufficiently long.  A modified version of the refined impedance transform adapted to the case of a zero mean wavelet is described below in \S\ref{sec-case}. The rate of convergence of (\ref{convergence}) depends on $\zeta$.  In cases where the convergence is especially fast one obtains an approximation $\imp_{w,\zeta_-}D(x_+)\cong\zeta(x_+)$. Extending the latter to the interior of the slab by replacing $x_+$ with $x<x_+$ results in the estimate
\begin{equation}\label{interior}
\zeta(x)\cong\imp_{w,\zeta_-}D(x).   
\end{equation}
Specializing to $W=\delta$, in which case $D=G_\zeta$, (\ref{interior}) yields an alternative to the classic approximation (\ref{standard}), namely
\begin{equation}\label{alternative}
\zeta(x)\cong\zeta_-\frac{1-\int_{-\infty}^{2x}G_\zeta}{1+\int_{-\infty}^{2x}G_{\zeta}}.  
\end{equation}
Both formulas involve the accumulation function of $G_\zeta$,
\begin{equation}\label{accumulation}
A(x)=\int_{-\infty}^{2x}G_\zeta=\int_0^{2x}G_\zeta.
\end{equation}
But remarkably, the simpler formula (\ref{alternative}) is fundamentally more accurate, as will be demonstrated in later sections.  Indeed, the two formulas agree to second order in the accumulation function: the expansion of (\ref{standard}) as a power series in $A$ is
\[
\zeta\cong\zeta_-e^{-2A}=\zeta_-\left(1-2A+2A^2-\textstyle\frac{4}{3}A^3+\cdots\right)
\]
while that of (\ref{alternative}) is
\[
\zeta\cong\zeta_-\frac{1-A}{1+A}=\zeta_-\left(1-2A+2A^2-2A^3+\cdots\right).
\]
Effectively, the classical approximation is agrees with (\ref{alternative}) in cases where $A$ is sufficiently small.  In many physically natural examples, however, this is not the case and the classical approximation fails (see \S\ref{sec-numerical} and \S\ref{sec-conclusion}). 
. 

Rearranging (\ref{alternative}) produces an estimate for the accumulation function of the Green's function in terms of $\zeta$,
\begin{equation}\label{accumulation-approximation}
A(x)\cong\frac{\zeta_--\zeta(x)}{\zeta_-+\zeta(x)}.  
\end{equation}
Based on this, we define the \emph{energy lag} of a given impedance profile $\zeta$ to be the function
\begin{equation}\label{energy-lag}
\sigma_\zeta(x)=\left|\frac{\zeta_--\zeta(x)}{\zeta_-+\zeta(x)}-\int_{0}^{2x}G_\zeta\right|.
\end{equation}
In physical terms, the energy lag measures the degree to which echoes that ultimately return to $x=0$ are delayed by reverberation (exclusively) within the interval $[0,x]$.  A small energy lag implies high accuracy of the refined impedance transform of the impulse response as an approximation to the impedance.  Theorem~\ref{thm-slab} is equivalent to the assertion that $\sigma_\zeta(x)\rightarrow0$ as $x\rightarrow\infty$.   

Differentiation of (\ref{accumulation-approximation}) yields a counterpart to the classical estimate (\ref{bad-approximation}), namely 
\begin{equation}\label{measurement-approximation}
G_\zeta(t)\cong-\zeta_-\frac{\zeta^\prime(t/2)}{\left(\zeta_-+\zeta(t/2)\right)^2}.
\end{equation}
Just as the refined impedance transform is more accurate than the standard estimate (\ref{standard}), the formula (\ref{measurement-approximation}) turns out to be substantially more accurate than the single scattering approximation (\ref{bad-approximation}), while also being very simple.   

\subsection{The case $w=0$\label{sec-case}}

Various procedures in seismic signal processing attempt to deconvolve the source wavelet $W$ from the data, thereby transforming $D$ to $G_\zeta$, at least approximately.  Formula (\ref{alternative}) applies to deconvolved data, as does Theorem~\ref{thm-slab} with $W=\delta$.  But in the absence of deconvolution the question of whether or not the wave form $W$ is zero mean becomes important.  In a typical experimental context the source wave form \emph{is} zero mean for the physical reason that this corresponds to no net displacement of the medium. 
%(net displacement being the integral of particle velocity).  
The following modified version of the refined impedance transform is tailored to the zero mean case.  

Let $k$ denote the least nonnegative integer such that $\int_{-\infty}^\infty s^kW(s)\,ds\neq0$. The assumption $\int W=0$ guarantees that $k\geq 1$.   Set 
\begin{equation}\label{V}
V(x)=-\int_{-\infty}^x\frac{(s-x)^{k-1}}{(k-1)!}W(s)\,ds.
\end{equation}
Then 
\begin{equation}\label{v}
v:=\int_{-\infty}^\infty V=\lim_{x\rightarrow\infty}\int_{-\infty}^x\frac{(s-x)^k}{k!}W(s)\,ds=\int_{-\infty}^\infty\frac{s^k}{k!}\,W(s)\,ds\neq0
\end{equation}
by definition of $k$.  Observe that for the measured data $D=\widetilde{W}\ast G_\zeta$, 
\begin{equation}\label{conv}
D^{(-k)}(t):=\int_{-\infty}^t\frac{(t-s)^{k-1}}{(k-1)!}D(s)\,ds=\widetilde{V}\ast G_\zeta(t).
\end{equation}
Thus $D^{(-k)}$ (the $k$-fold antiderivative of $D$) has the structure of measured data corresponding to (virtual) source wave form $V$, which is non zero mean by construction.  One can then apply the refined impedance transform $\imp_{v,\zeta_-}$ to $D^{(-k)}$ to obtain an estimate for the impedance function $\zeta$, extending Theorem~\ref{thm-slab} as follows.  
\begin{thm}\label{thm-general}
Let $D$ denote the reflection data for (\ref{wave}) as defined by (\ref{data}), where the source wave form $W$ is not identically zero, and suppose that $\int W=0$.  Define $v$ and $D^{(-k)}$ according to (\ref{v}) and (\ref{conv}), respectively, so that $k\geq1$ denotes the least nonzero moment of $W$. Then 
\begin{equation}\label{general-convergence}
\imp_{v,\zeta_-}D^{(-k)}(x)\rightarrow\zeta_+\qquad\mbox{ as }\qquad x\rightarrow\infty.
\end{equation}
\end{thm}
The analogous estimate to (\ref{interior}) in the zero mean case is 
\begin{equation}\label{analogous}
\zeta(x)\cong\imp_{v,\zeta_-}D^{(-k)}(x).
\end{equation}

\section{Derivation and proof\label{sec-derivation}}
Theorem~\ref{thm-slab} rests on two key ideas: approximation of a given impedance function $\zeta$ by step functions, and an explicit formula for the integral of the reflection Green's function of a step function.  More precisely, among those impedance functions which take constant values to the left and right of an interval $[x_-,x_+]$ for some $0<x_-<x_+$, and which are bounded and bounded away from zero, we consider those that are uniform limits of step functions, i.e., \emph{regulated} functions.  Regulated functions are a broad class that includes, for example, all functions of bounded variation (see \cite[Ch.~VII]{Di:1960}).  In particular, regulated functions may have discontinuities---possibly infinitely many. %(Regulated functions cannot have the type of singularity that $\sin(1/x)$ has at zero, however.)  
Thus, from the point of view of modelling a physical medium, $\zeta$ can be essentially arbitrary.  

A crucial fact from the general theory of time dependent linear equations is that if $\zeta^P\rightarrow\zeta$ uniformly for some sequence of partitions $P$, then $\int G_{\zeta^P}\rightarrow\int G_\zeta$ (for details see \cite[Lemma~3.1]{KiRi:2014} or in a more general context \cite[Ch.~8]{LiMa:1972}).  This ensures that results established for step functions $\zeta$ carry over to the more general class of regulated functions.  

We recall the key formula for $\int G_\zeta$ in the case where $\zeta$ is a step function.   Let $\zeta(x\pm)$ denote one sided limits $\lim_{y\rightarrow x\pm}\zeta(y)$.
If $\zeta$ is a step function, the function
\[
\frac{\zeta(x-)-\zeta(x+)}{\zeta(x-)+\zeta(x+)}\qquad(x\in\real),
\]
is non-zero only at jump points, which may be listed in order as $x_1<\cdots<x_n$; assume without loss of generality that $x_1>x_0=0$. Writing 
\begin{equation}\label{reflectivities}
r_j=\frac{\zeta(x_j-)-\zeta(x_j+)}{\zeta(x_j-)+\zeta(x_j+)}\quad\mbox{ and }\quad \tau_j=2(x_j-x_{j-1})\qquad(1\leq j\leq n),
\end{equation}
and setting $\tau=(\tau_1,\ldots,\tau_n)$ and $r=(r_1,\ldots,r_n)$, the pair $(\tau,r)$ encodes $\zeta$.  In this case the reflection Green's function for (\ref{wave}) has the form 
\begin{equation}\label{form}
G_\zeta(t)=G^{(\tau,r)}(t)=\sum_{k\in\{1\}\times\integer^{n-1}}\mathfrak{a}(r,k)\delta(t-\langle\tau,k\rangle),
\end{equation}
where $\mathfrak{a}(r,k)$ are explicitly given polynomials in $r$ \cite[Thm.~4]{Gi:SIAP2014}.  Moreover, \cite[Thm.~2]{Gi:IPI2017} ensures that, provided each $|r_j|<1$, the amplitudes $\mathfrak{a}(r,k)$ are absolutely summable, and 
\begin{equation}\label{sum}
\sum_{k\in\{1\}\times\integer^{n-1}}\mathfrak{a}(r,k)=\tanh\left(\sum_{j=1}^n\tanh^{-1}r_j\right).
\end{equation}

\subsection{Proof of Theorem~\ref{thm-slab}}

Let $\zeta$ be a regulated function with $\alpha$ and $\beta$ constants such that $0<\alpha<\zeta(x)<\beta$ $(x\in\real)$.  Suppose first that $\zeta$ is a step function.  Then the reflection Green's function may be expressed in the form $G_\zeta(t)=\sum_{i=1}^\infty a_i\delta(t-t_i)$, with each $a_i\neq0$ and $t_i<t_{i+1}$.  Defining reflectivities $r_j$ as above in (\ref{reflectivities}) 
we then have that $|r_j|<(\beta-\alpha)/(\beta+\alpha)<1$ $(1\leq j\leq n)$. Therefore (\ref{sum}) implies that
\[
\sum_{i=1}^\infty a_i=\tanh\left(\sum_{j=1}^n\tanh^{-1}r_j\right),
\]
with the series on the left converging absolutely.  The right-hand sum telescopes since
\[
r_j=\frac{\zeta(x_j-)-\zeta(x_j+)}{\zeta(x_j-)+\zeta(x_j+)}=\tanh\left(\frac{1}{2}\left(\log\zeta(x_j-)-\log\zeta(x_j+)\right)\right)
\]
to yield
\begin{equation}\label{zeta-sum}
\int_{-\infty}^\infty G_\zeta=\sum_{i=1}^\infty a_i=\tanh\left(\frac{1}{2}\left(\log\zeta_--\log\zeta_+\right)\right)=\frac{\zeta_--\zeta_+}{\zeta_-+\zeta_+}.
\end{equation}

Now consider the general case of a regulated function $\zeta$.  
Since $\zeta$ is regulated, there exists a sequence of partitions $P$ such that $\zeta^P\rightarrow\zeta$ uniformly, and for every $P$, $\zeta^P_-=\zeta_-$ and $\zeta^P_+=\zeta_+$.  It follows by (\ref{zeta-sum}) that 
\[
\int_{-\infty}^\infty G_{\zeta^P}=\frac{\zeta_--\zeta_+}{\zeta_-+\zeta_+}
\]
is independent of $P$.  Therefore   
\[
\int_{-\infty}^\infty G_\zeta=\lim_P\int_{-\infty}^\infty G_{\zeta^P}=\frac{\zeta_--\zeta_+}{\zeta_-+\zeta_+}. 
\]
Rearranging terms,
\[
\zeta_+=\zeta_-\frac{1-\int_{-\infty}^\infty G_\zeta}{1+\int_{-\infty}^\infty G_\zeta}=\zeta_-\frac{\left(\int_{-\infty}^\infty W\right)-\left(\int_{-\infty}^\infty W\right)\int_{-\infty}^\infty G_\zeta}{\left(\int_{-\infty}^\infty W\right)+\left(\int_{-\infty}^\infty W\right)\int_{-\infty}^\infty G_\zeta}=\zeta_-\frac{w-\int_{-\infty}^\infty D}{w+\int_{-\infty}^\infty D}=\lim_{x\rightarrow\infty}\imp_{w,\zeta_-}D,
\]
completing the proof of Theorem~\ref{thm-slab}.  

Before giving some numerical examples, we note the behaviour of the refined impedance transform with respect to scaling and dilation of the impedance function, and corresponding dilation of the source wave form.  
\begin{prop}\label{prop-dilation}
Let $a,b>0$ and suppose $\zeta^1,\zeta^2,W^1,W^2$ are such that 
\[
\zeta^2(x)=a\zeta^1(bx)\mbox{ and } W^2(x)=bW^1(bx)\qquad(x\in\real).
\]
Write $w=\int_{-\infty}^\infty W^1=\int_{-\infty}^\infty W^2$ and set $\xi^j=\imp_{w,\zeta_-^j}\widetilde{W^j}\ast G_{\zeta^j}$ for $j=1,2$.  Then 
\[
\xi^2(x)=a\xi^1(bx)\qquad(x\in\real).
\]
\end{prop}
\begin{pf}
Suppose first that $\zeta^1$ is a step function. Then, as above, the function
\[
\frac{\zeta^1(x-)-\zeta^1(x+)}{\zeta^1(x-)+\zeta^1(x+)},
\]
is non-zero only at jump points, which may be listed in order as $x^1_1<\cdots<x^1_n$, with $0<x^1_-\leq x^1_1<x^1_n\leq x^1_+$.  Write $x^1_0=0$ and set 
\[
r^1_j=\frac{\zeta(x^1_j-)-\zeta(x^1_j+)}{\zeta(x^1_j-)+\zeta(x^1_j+)}\quad\mbox{ and }\quad \tau^1_j=2(x^1_j-x^1_{j-1})\qquad(1\leq j\leq n).
\]
Writing $r^1=(r^1_1,\ldots,r^1_n)$ and $\tau^1=(\tau^1_1,\ldots,\tau^1_n)$, we have $G_{\zeta^1}=G^{(\tau^1,r^1)}$.  The same definitions with $\zeta^2$ in place of $\zeta^1$ yield $G_{\zeta^2}=G^{(\tau^2,r^2)}$, where $r^2=r^1$ and $\tau^2=\frac{1}{b}\tau^1$.  
It follows from (\ref{form}) that 
\[
\begin{split}
G_{\zeta^2}(t)&=\sum_{k\in\{1\}\times\integer_+^{n-1}}\mathfrak{a}(r^2,k)\delta\left(t-\langle\tau^2,k\rangle\right)\\
&=\sum_{k\in\{1\}\times\integer_+^{n-1}}\mathfrak{a}(r^1,k)\delta\left(t-\frac{1}{b}\langle\tau^1,k\rangle\right)\\
&=\sum_{k\in\{1\}\times\integer_+^{n-1}}\mathfrak{a}(r^1,k)b\delta\left(bt-\langle\tau^1,k\rangle\right)\\
&=bG_{\zeta^1}(bt).
\end{split}
\]
(Note that the equation $\sigma(x)=\delta(x/b)$, properly interpreted in terms of test functions, defines a distribution $\sigma=b\delta$.)  Thus,
\[
\int_{-\infty}^{2x}\widetilde{W^2}\ast G_{\zeta^2}=\int_{-\infty}^{2x}\widetilde{bW^1(b\,\cdot)}\ast bG_{\zeta^1}(b\,\cdot)=\int_{-\infty}^{2bx}\widetilde{W^1}\ast G_{\zeta^1},
\]
which, given that $\zeta^2_-=a\zeta^1_-$, yields the desired relation $\xi^2(x)=a\xi^1(bx)$ $(x\in\real)$.  

The result for general regulated $\zeta^1$ then follows using approximation by step functions.  
\end{pf}

\section{Numerical examples\label{sec-numerical}}
Standard finite different schemes cannot guarantee accuracy of numerical computation of the reflection Green's function for (\ref{wave}) since they require smooth initial data. The present paper exploits a novel method that guarantees accuracy even for a purely distributional source wave form and discontinuous impedance---see \cite{Gi:NMPDE2018} for details.

\subsection{The cases $W=\delta$ and Gaussian $W$}

\begin{figure}[p]
\fbox{
\parbox{221pt}{
\includegraphics[clip,trim=0in 0in 0in 0in, width=1.5in]{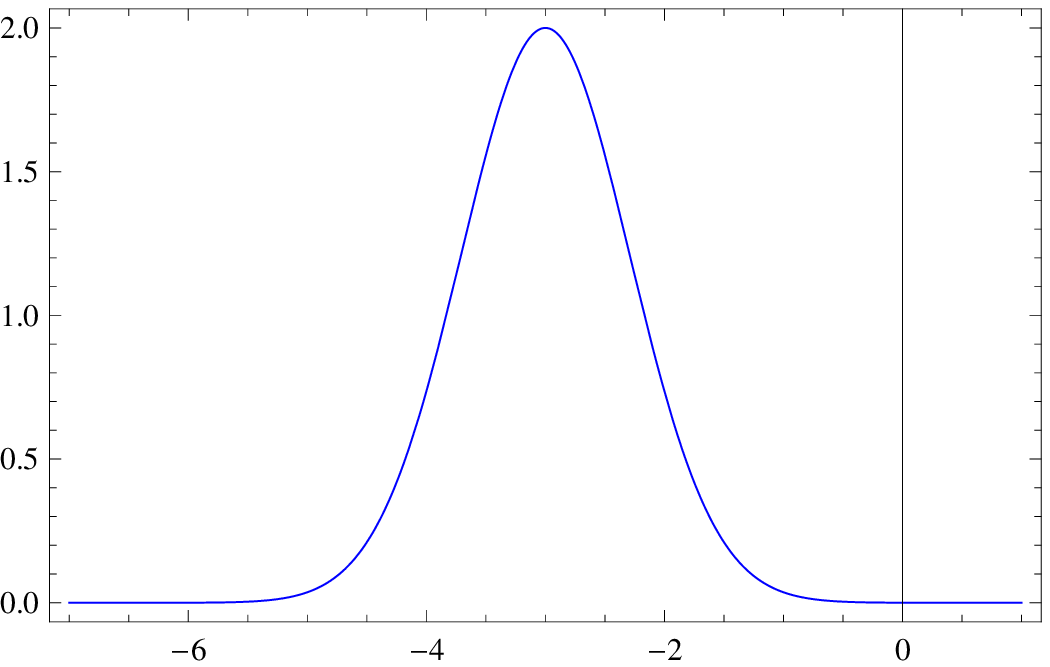}\hspace*{5pt}\includegraphics[clip,trim=0in 0in 0in 0in, width=1.5in]{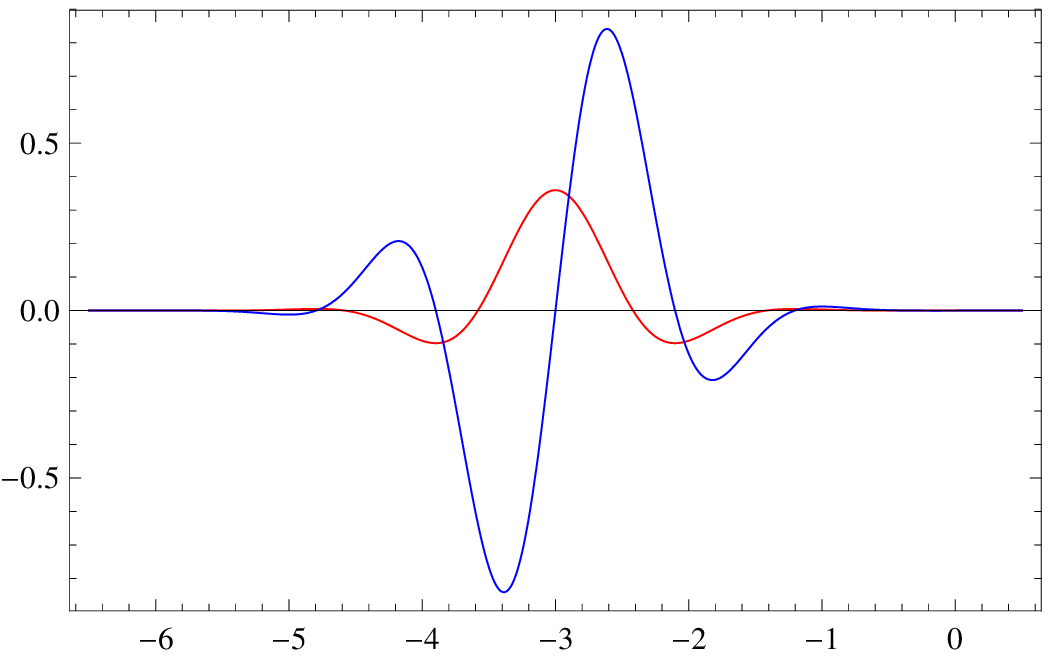}
\caption{Two source wave forms.  Left: A Gaussian wave.  Right: A zero mean source waveform (blue), plotted together with the negative of its antiderivative (red), which serves as the virtual wave form in the modified transform.}\label{fig-wave-form}
}
}
\end{figure}

The present section compares the modified transform to the standard approximation applied to reflection data coming from four different impedance profiles.  First the transforms are applied to the reflection Green's function, and then to undeconvolved data corresponding to a Gaussian source wave form depicted in Figure~\ref{fig-wave-form}.

\begin{figure}[p]
\fbox{
\parbox{452pt}{
\includegraphics[clip,trim=0in 0in 0in 0in, width=2in]{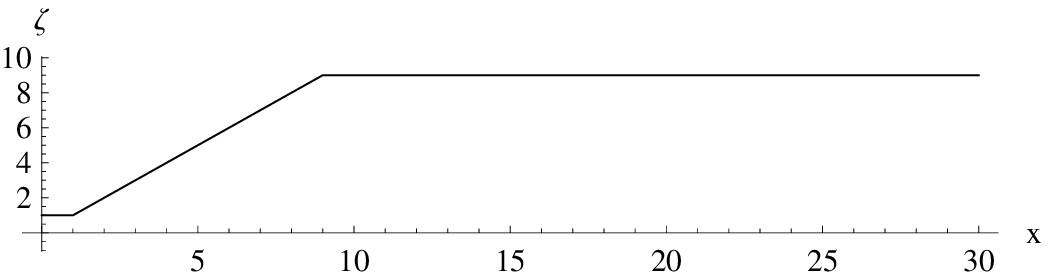}\hspace*{10pt}\includegraphics[clip,trim=0in -10pt 0in 0in, width=2in]{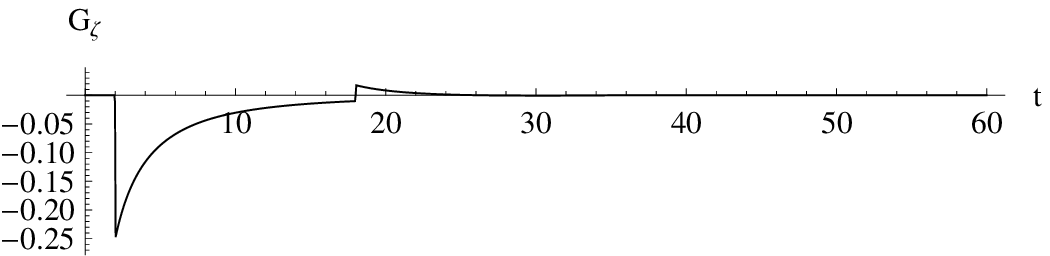}\hspace*{10pt}\includegraphics[clip,trim=0in -10pt 0in 0in, width=2in]{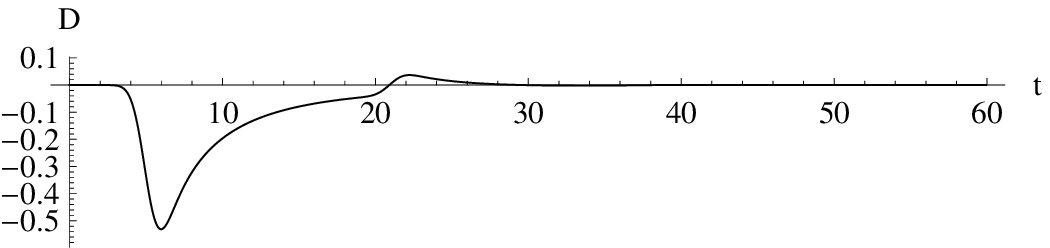}
\includegraphics[clip,trim=0in 0in 0in 0in, width=2in]{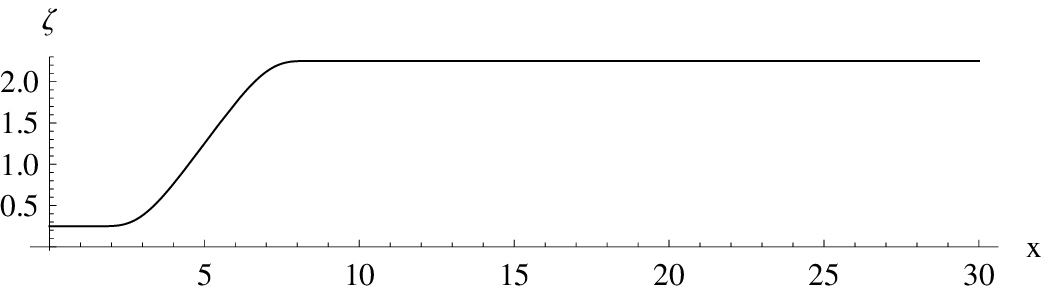}\hspace*{10pt}\includegraphics[clip,trim=0in -10pt 0in 0in, width=2in]{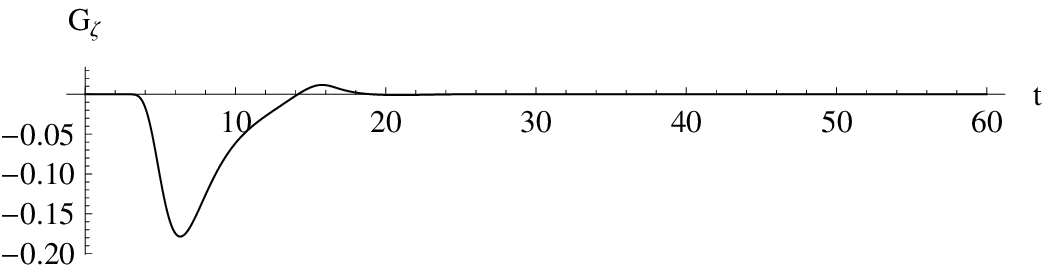}\hspace*{10pt}\includegraphics[clip,trim=0in -10pt 0in 0in, width=2in]{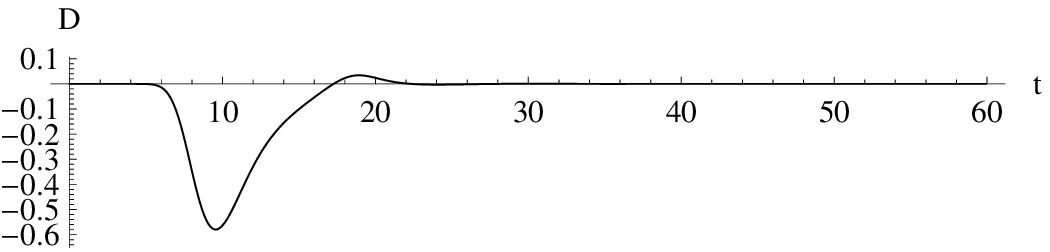}
\includegraphics[clip,trim=0in 0in 0in 0in, width=2in]{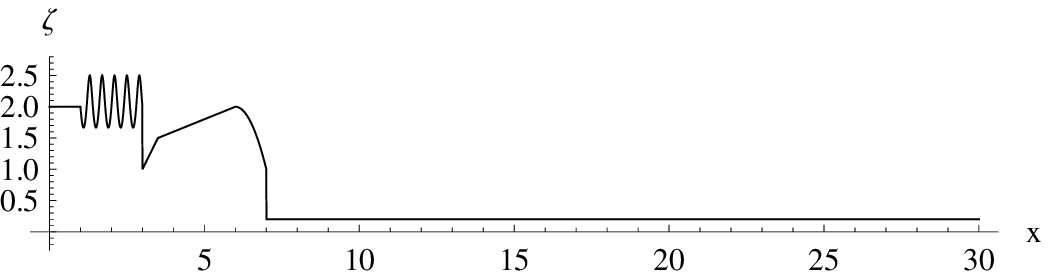}\hspace*{10pt}\includegraphics[clip,trim=0in -10pt 0in 0in, width=2in]{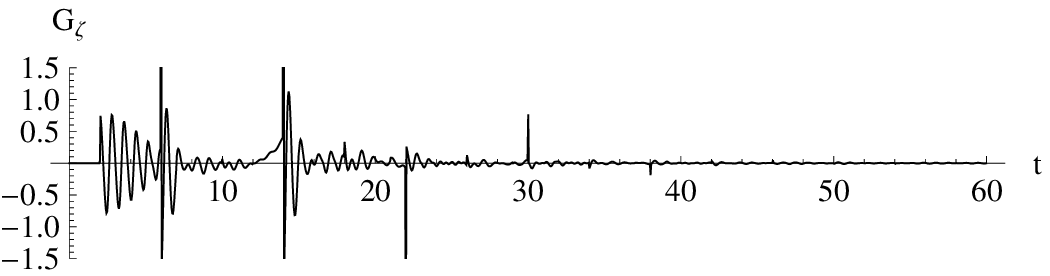}\hspace*{10pt}\includegraphics[clip,trim=0in -10pt 0in 0in, width=2in]{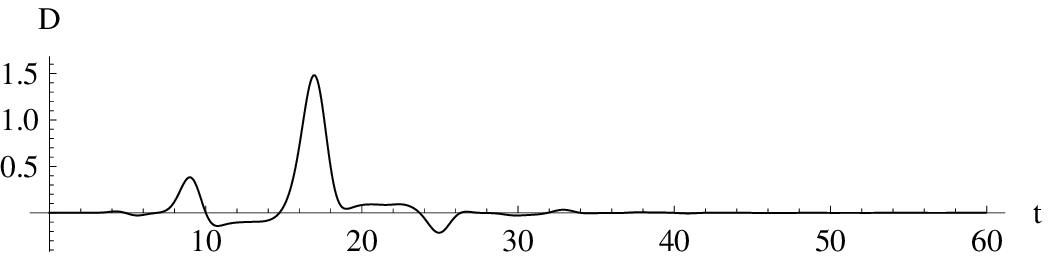}
\includegraphics[clip,trim=0in 0in 0in 0in, width=2in]{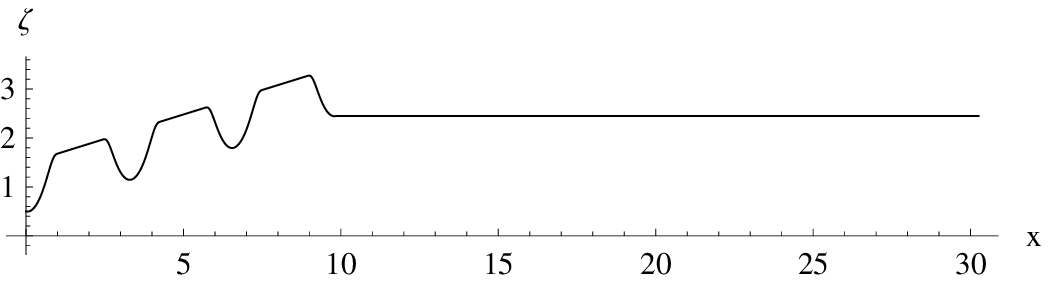}\hspace*{10pt}\includegraphics[clip,trim=0in -10pt 0in 0in, width=2in]{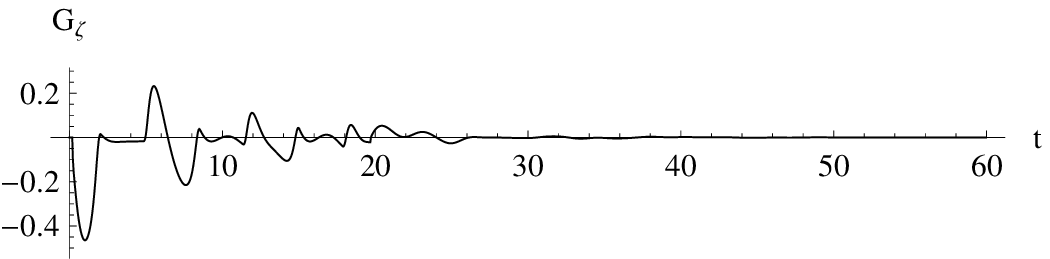}\hspace*{10pt}\includegraphics[clip,trim=0in -10pt 0in 0in, width=2in]{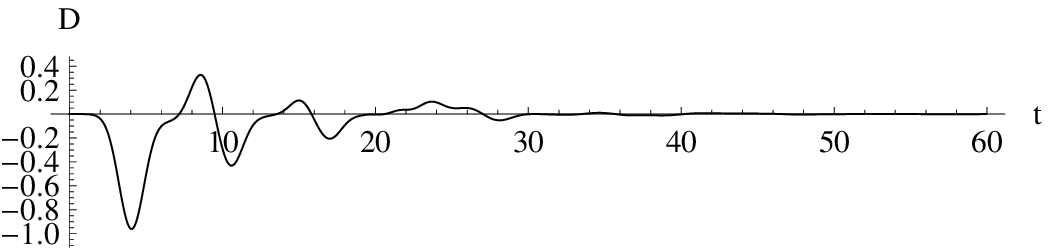}
\caption{Impedance profiles (on the left) with their reflection Green's functions (middle) and undeconvolved reflection data (right) corresponding to a Gaussian source wave.}\label{fig-profiles}
}
}
\end{figure}

\begin{figure}[p]
\fbox{
\parbox{226pt}{
\includegraphics[clip,trim=0in 0in 0in 0in, width=1.5in]{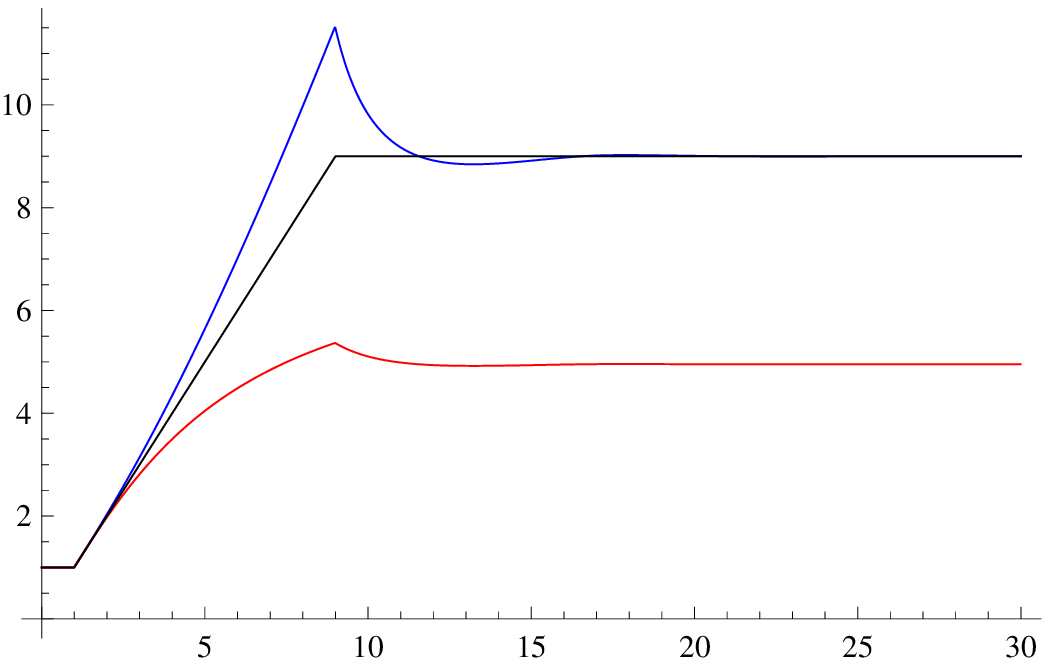}\hspace*{10pt}\includegraphics[clip,trim=0in 0in 0in 0in, width=1.5in]{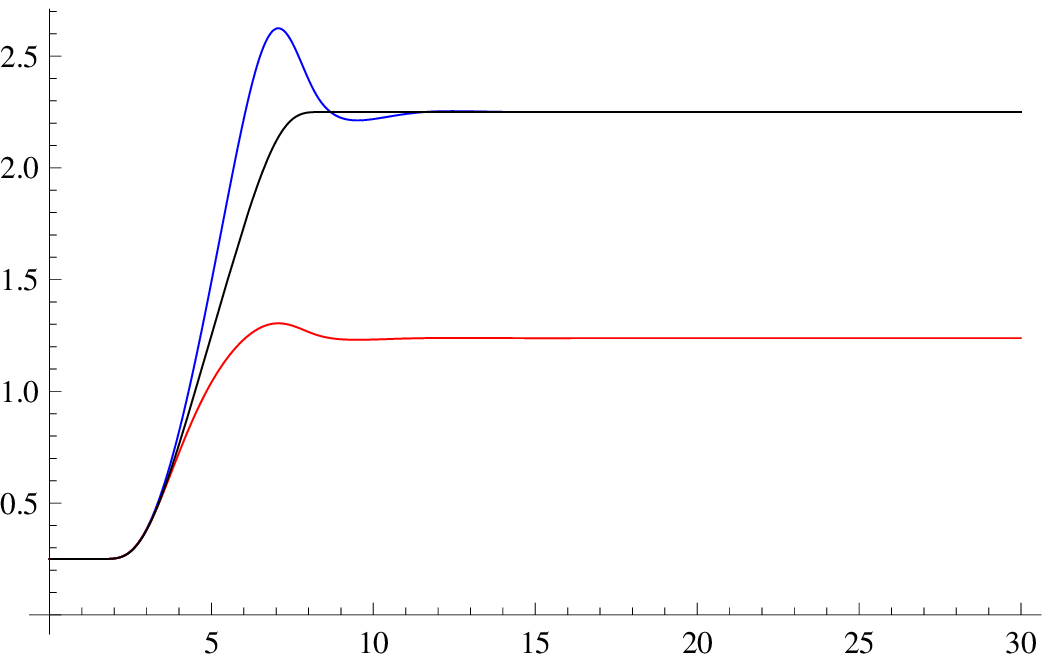}\\
\includegraphics[clip,trim=0in 0in 0in 0in, width=1.5in]{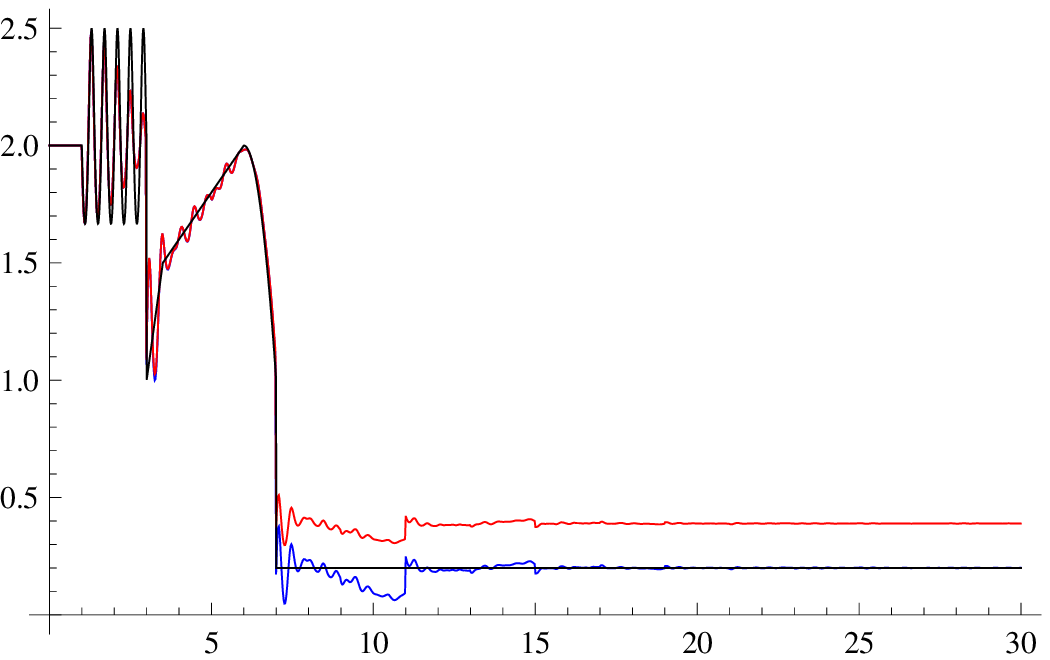}\hspace*{10pt}\includegraphics[clip,trim=0in 0in 0in 0in, width=1.5in]{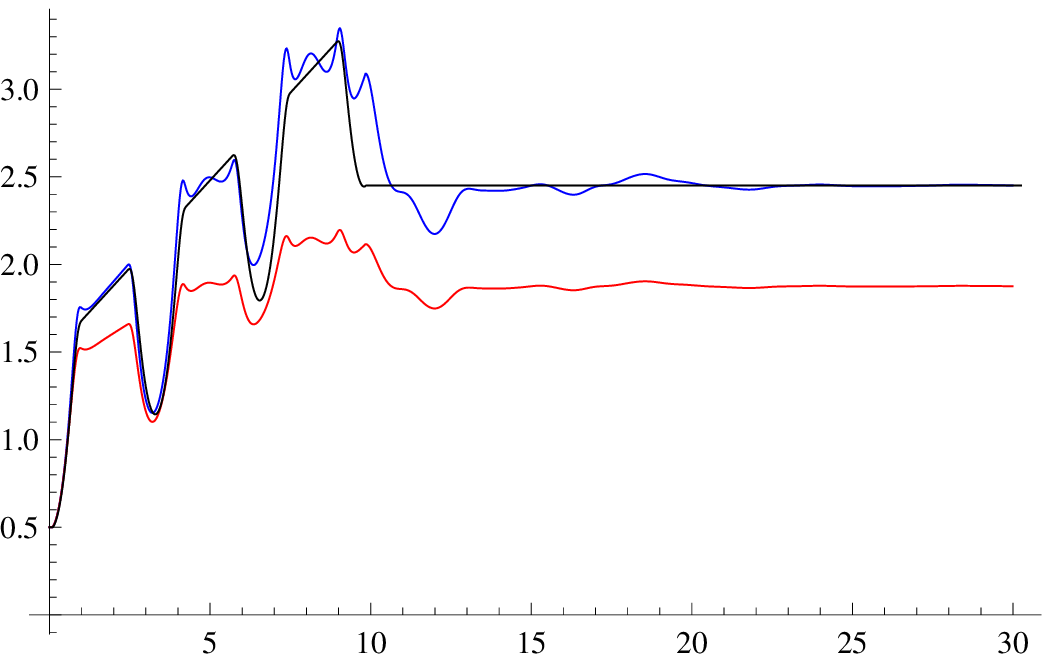}
\caption{Comparison of the modified transform (blue curve) to the standard approximation (red) applied to the reflection Green's function for each of the four impedance profiles in Figure~\ref{fig-profiles}.  The original impedance profile (black) is included for comparison.  In each case the modified transform eventually converges to the value $\zeta_+$ in accordance with Theorem~\ref{thm-slab}, whereas the standard approximation does not.}\label{fig-greens}
}
}
%\end{figure}
%
%\begin{figure}[p]
\fbox{
\parbox{215pt}{
\includegraphics[clip,trim=0in 0in 0in 0in, width=1.45in]{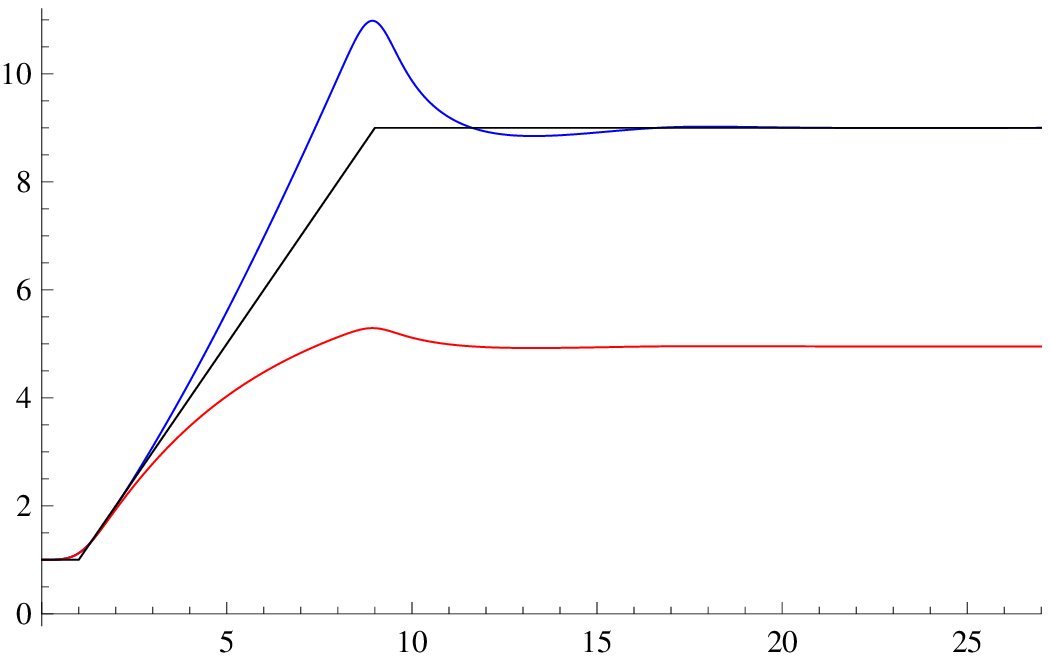}\hspace*{10pt}\includegraphics[clip,trim=0in 0in 0in 0in, width=1.45in]{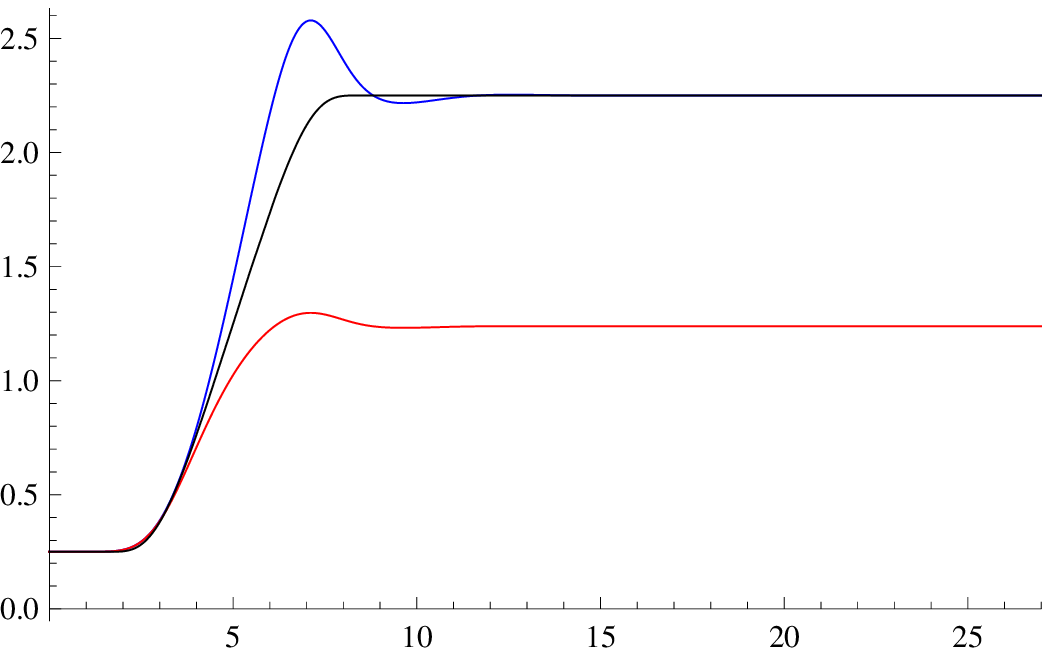}\\
\includegraphics[clip,trim=0in 0in 0in 0in, width=1.45in]{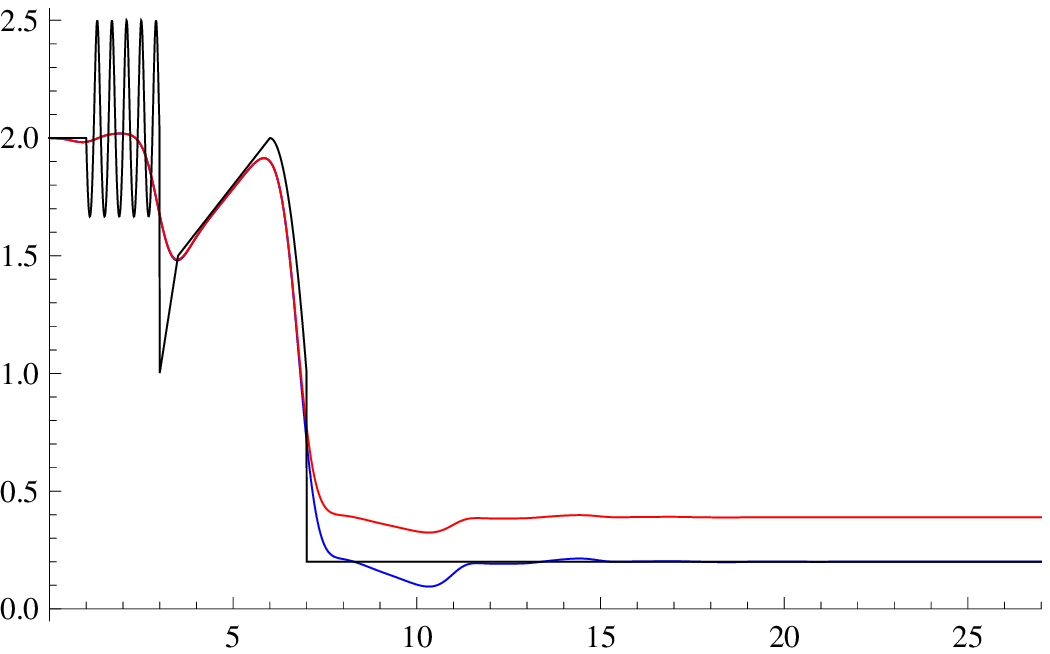}\hspace*{10pt}\includegraphics[clip,trim=0in 0in 0in 0in, width=1.5in]{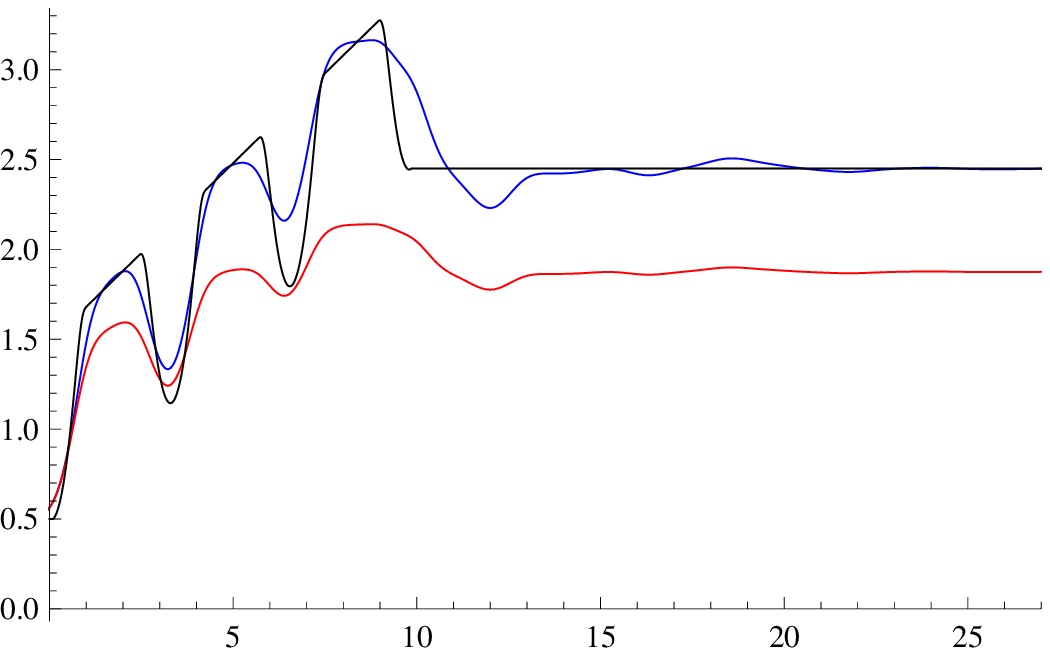}
\caption{The case of a Gaussian source wave form.  Comparison of the modified transform (blue curve) to the standard approximation (red) applied to the undeconvolved data for each of the four impedance profiles in Figure~\ref{fig-profiles}, with the original impedance profile (black) for reference.  As with the reflection Green's function, the modified transform eventually converges to the value $\zeta_+$.}\label{fig-wave}
}
}
\end{figure}

The impedance profiles, the reflection Green's functions and the undeconvolved data are depicted in Figure~\ref{fig-profiles}.  Figure~\ref{fig-greens} displays the approximate reconstructions of the impedance profiles obtained by applying the modified transform and standard approximation (\ref{standard}) to the reflection Green's functions.  Figure~\ref{fig-wave} shows the corresponding results for undeconvolved data with Gaussian source. 

\setlength{\fboxsep}{2.5pt}

\subsection{The case $w=0$}

\begin{figure}[p]
\fbox{
\parbox{447pt}{
\includegraphics[clip,trim=0in 0in 0in 0in, width=1.5in]{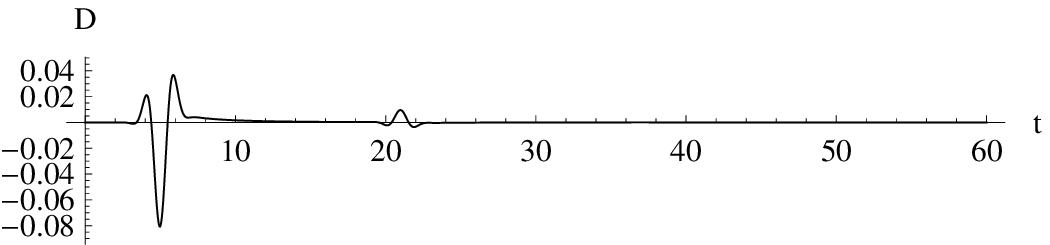}\hspace*{5pt}\includegraphics[clip,trim=0in 0in 0in 0in, width=1.5in]{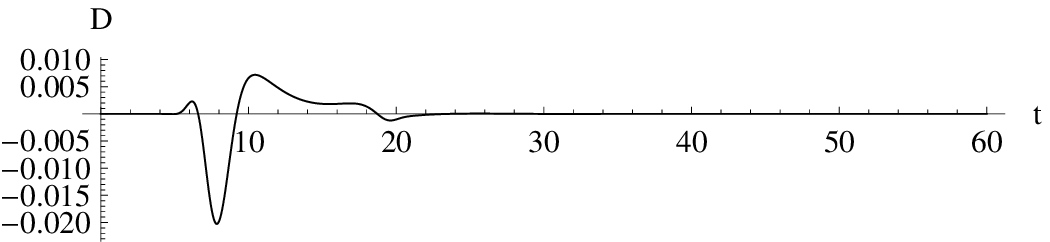}\hspace*{5pt}\includegraphics[clip,trim=0in 0in 0in 0in, width=1.5in]{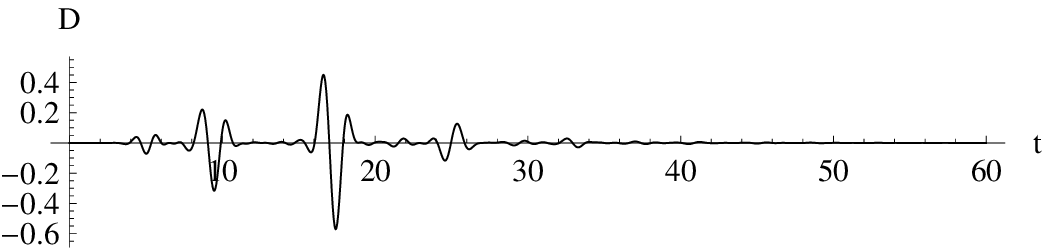}\hspace*{5pt}\includegraphics[clip,trim=0in 0in 0in 0in, width=1.5in]{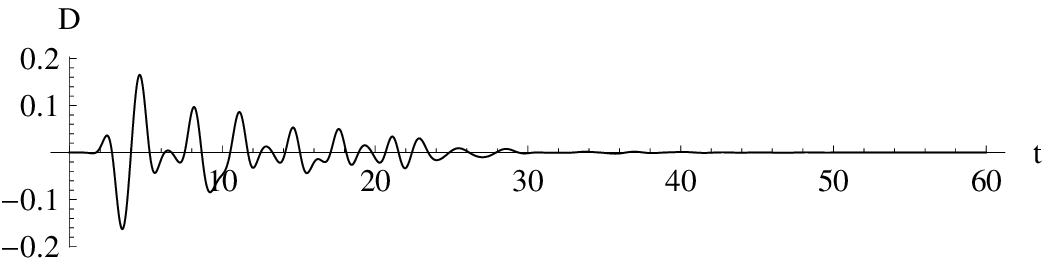}\\[10pt] 
\includegraphics[clip,trim=0in 0pt 0in 0in, width=1.5in]{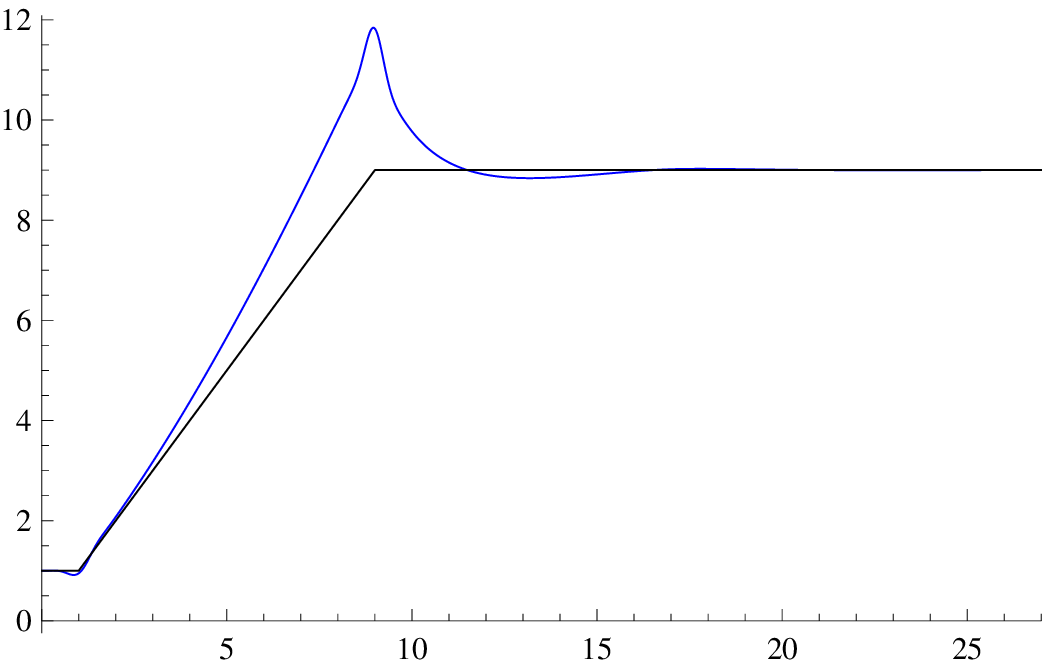}\hspace*{5pt}\includegraphics[clip,trim=0in 0pt 0in 0in, width=1.5in]{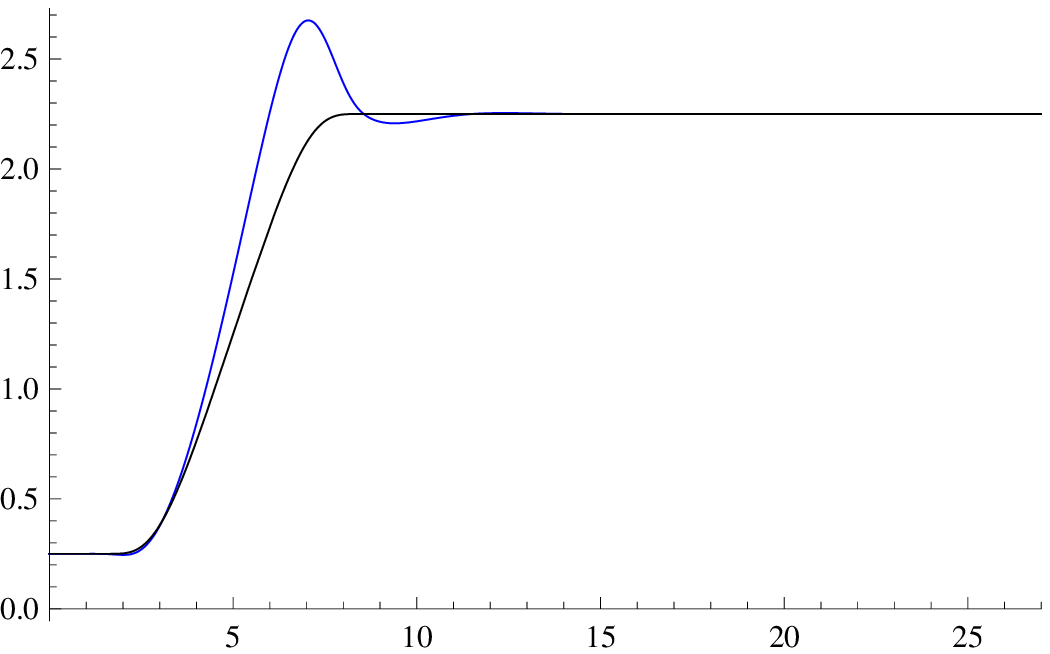}\hspace*{5pt}\includegraphics[clip,trim=0in 0pt 0in 0in, width=1.5in]{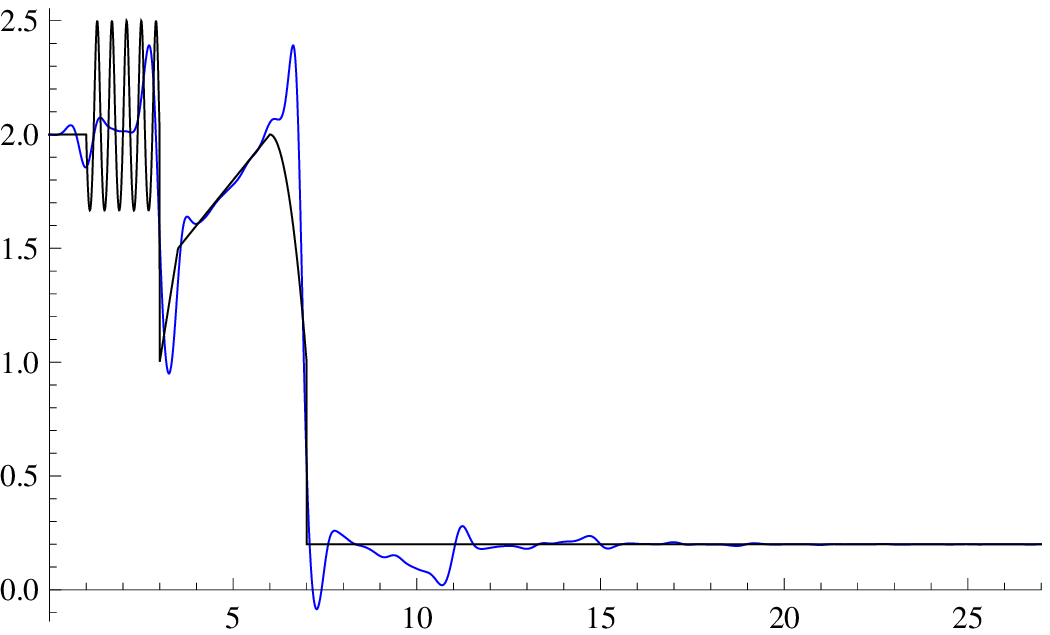}\hspace*{5pt}
\includegraphics[clip,trim=0in 0pt 0in 0in, width=1.5in]{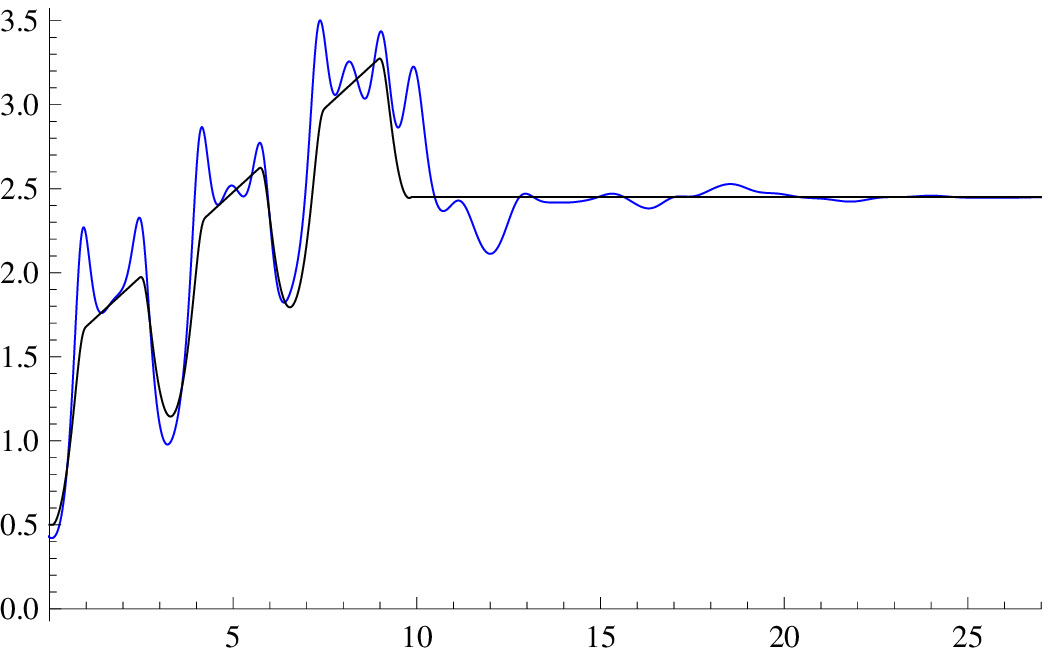}\caption{Impedance profiles (blue, bottom row) reconstructed from zero mean source wave form reflection data (top row) using the modified transform.}\label{fig-zero}
}
}
\end{figure}

The right-hand, blue wave form depicted in Figure~\ref{fig-wave-form} has zero mean and non-zero first moment.  The corresponding data and reconstructions obtained using the modified transform (with $k=1$) are depicted in Figure~\ref{fig-zero}.  These reconstructions are comparable in quality and accuracy to those obtained with a Gaussian source.  In other words, the modified transform effectively compensates for the zero mean source wave.   

\subsection{Noisy data\label{sec-noisy}}

Figures~\ref{fig-noisy-gaussian} and \ref{fig-noisy-zero-mean} below illustrate the effects of noisy data on the modified transform.  In the first case, 10\% noise is added to undeconvolved data for which the source wave form is non zero mean.  In the second case, 5\% noise is added to undeconvolved data for which the source wave form has zero mean.  
\begin{figure}[p]
\fbox{
\parbox{447pt}{
\includegraphics[clip,trim=0in 0in 0in 0in, width=1.5in]{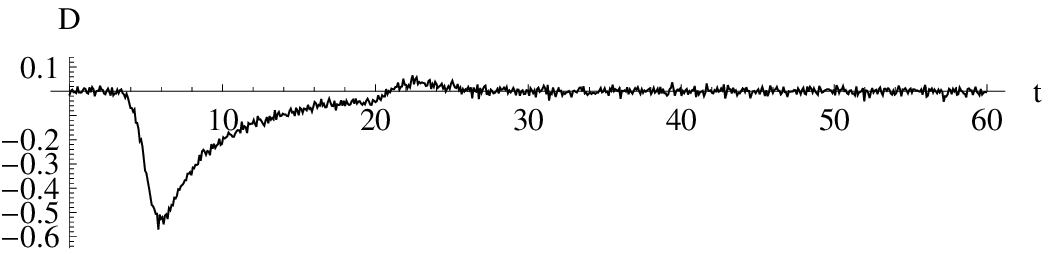}\hspace*{5pt}\includegraphics[clip,trim=0in 0in 0in 0in, width=1.5in]{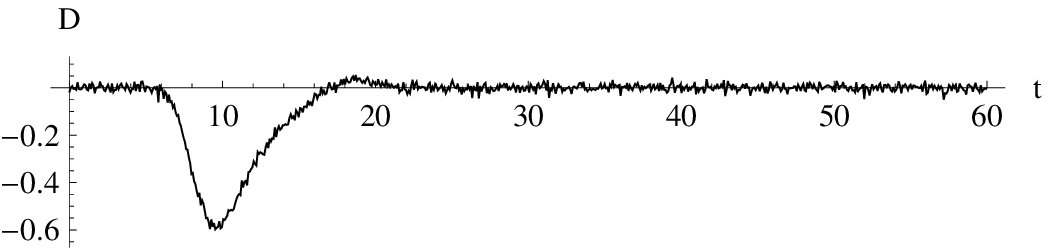}\hspace*{5pt}\includegraphics[clip,trim=0in 0in 0in 0in, width=1.5in]{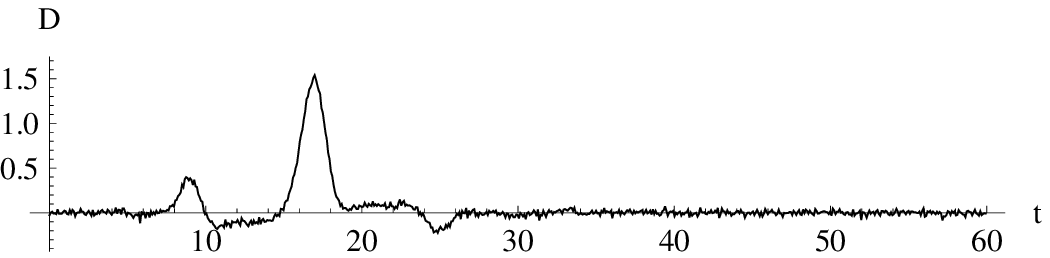}\hspace*{5pt}\includegraphics[clip,trim=0in 0in 0in 0in, width=1.5in]{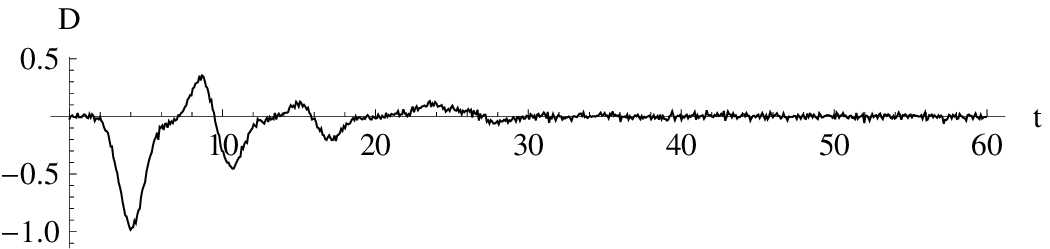}\\[10pt] 
\includegraphics[clip,trim=0in 0pt 0in 0in, width=1.5in]{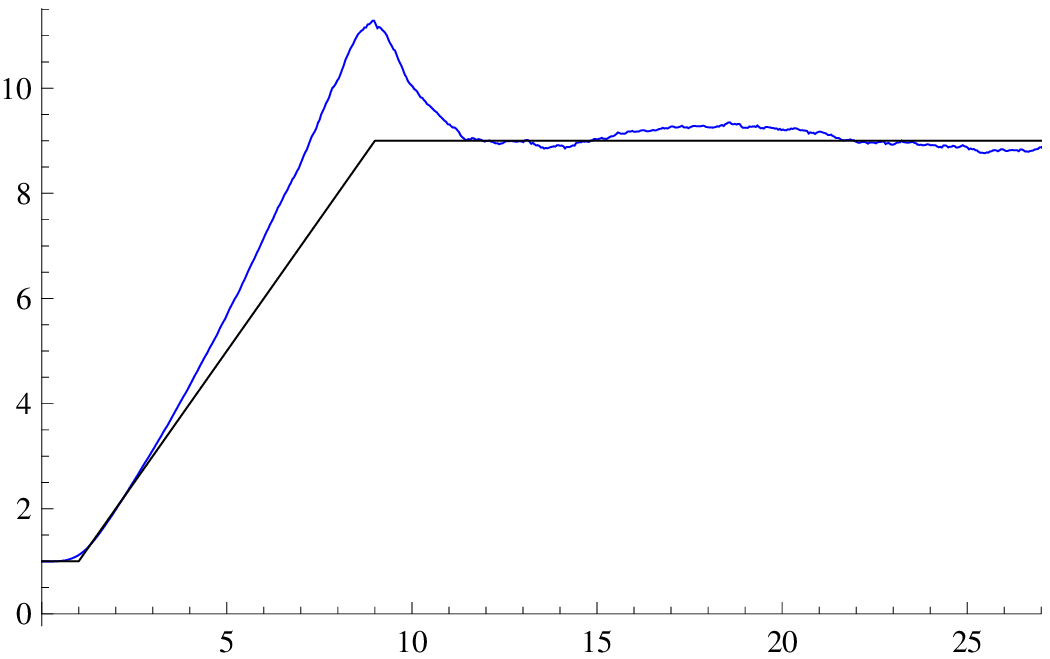}\hspace*{5pt}\includegraphics[clip,trim=0in 0pt 0in 0in, width=1.5in]{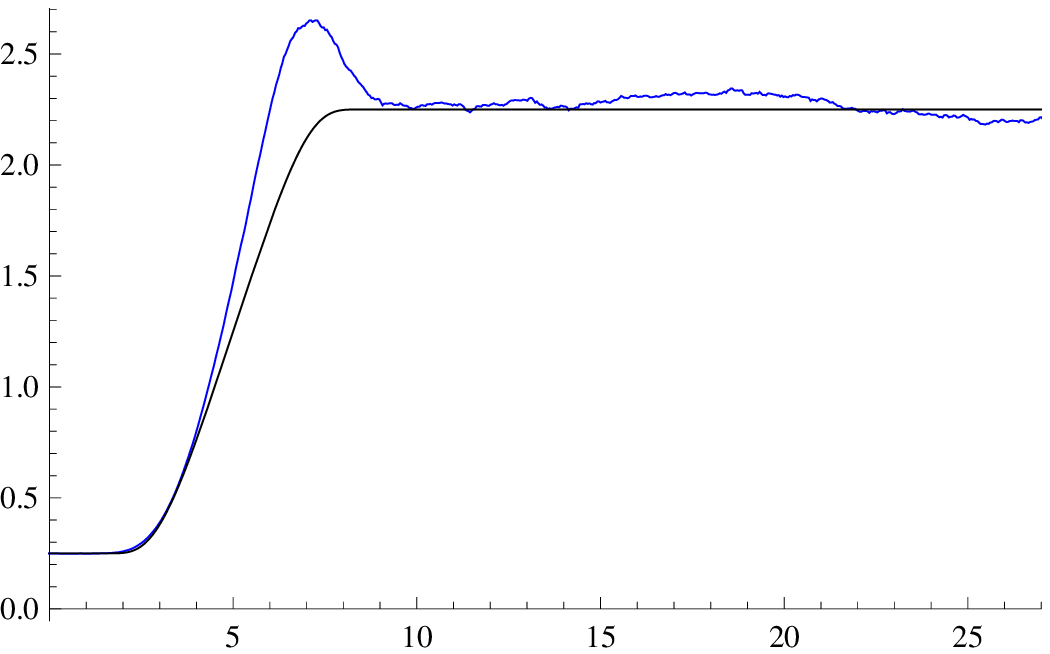}\hspace*{5pt}\includegraphics[clip,trim=0in 0pt 0in 0in, width=1.5in]{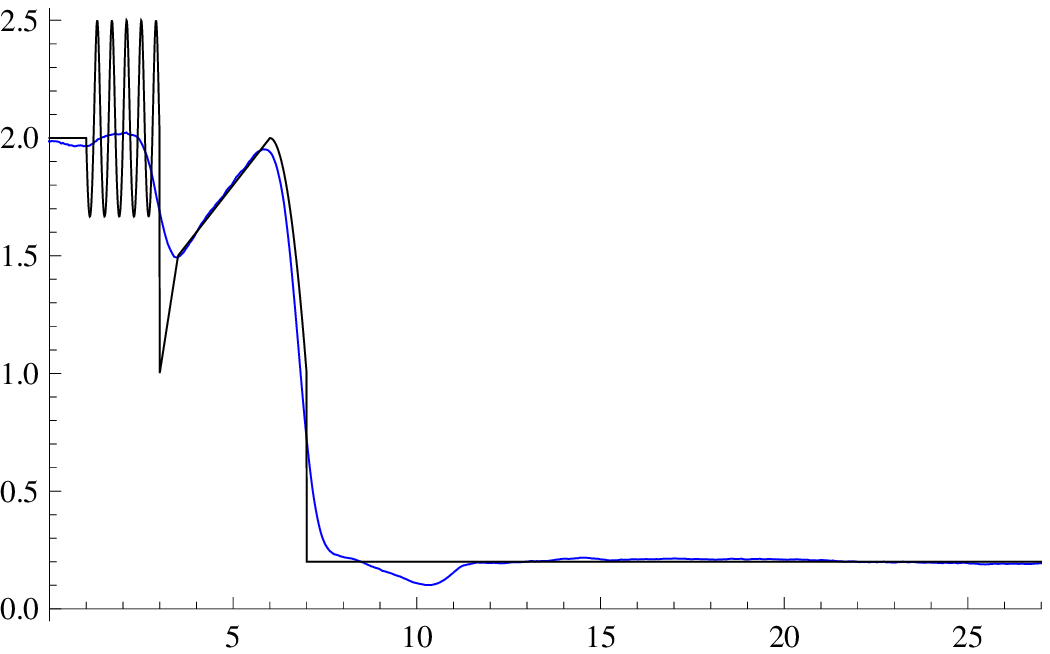}\hspace*{5pt}
\includegraphics[clip,trim=0in 0pt 0in 0in, width=1.5in]{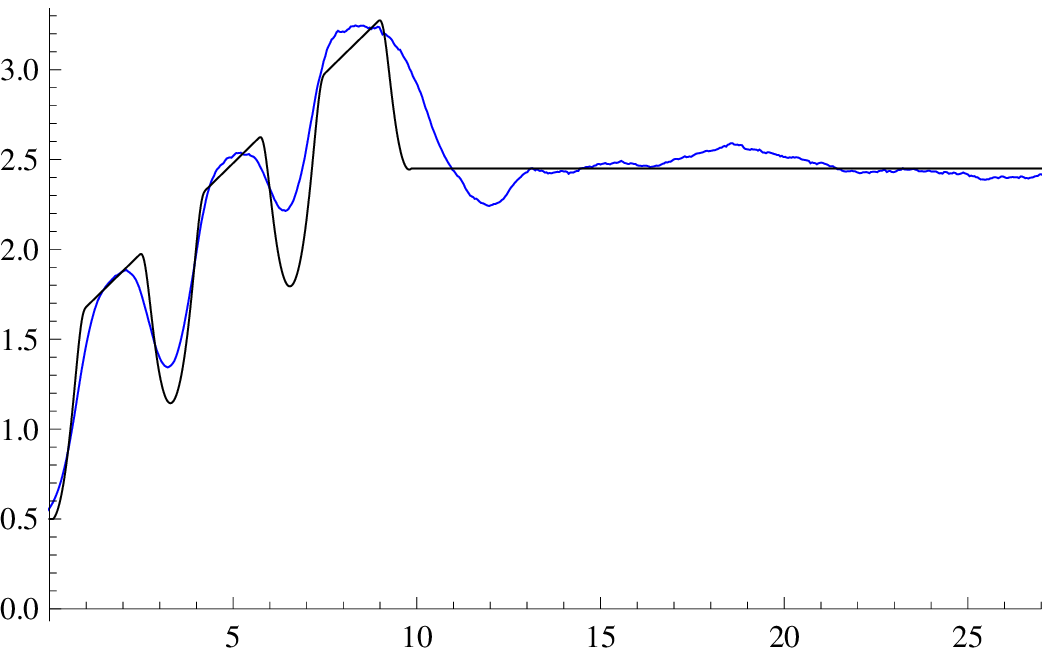}\caption{The effect of 10\% i.i.d.~Gaussian noise. Impedance profiles (blue, bottom row) reconstructed from noisy reflection data (top row) using the modified transform. The source wave is plotted on the left of Figure~\ref{fig-wave-form}. }\label{fig-noisy-gaussian}
}
}
\end{figure}

\begin{figure}[p]
\fbox{
\parbox{447pt}{
\includegraphics[clip,trim=0in 0in 0in 0in, width=1.5in]{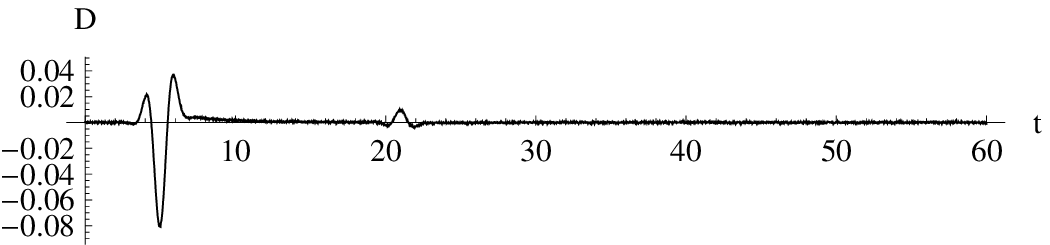}\hspace*{5pt}\includegraphics[clip,trim=0in 0in 0in 0in, width=1.5in]{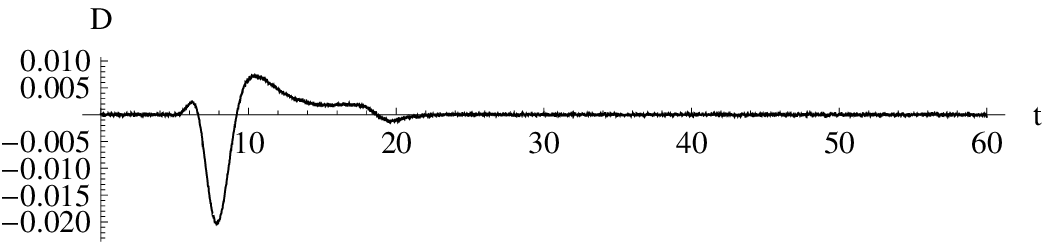}\hspace*{5pt}\includegraphics[clip,trim=0in 0in 0in 0in, width=1.5in]{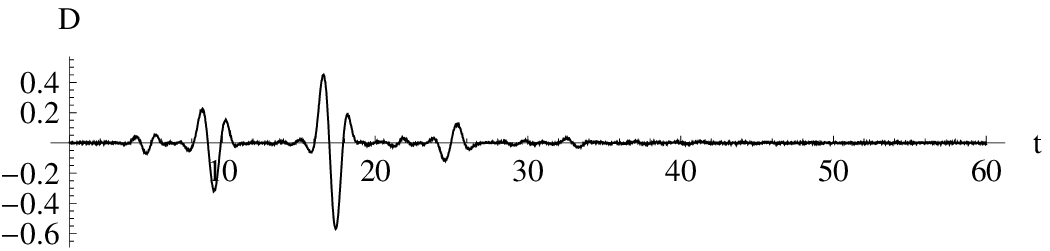}\hspace*{5pt}\includegraphics[clip,trim=0in 0in 0in 0in, width=1.5in]{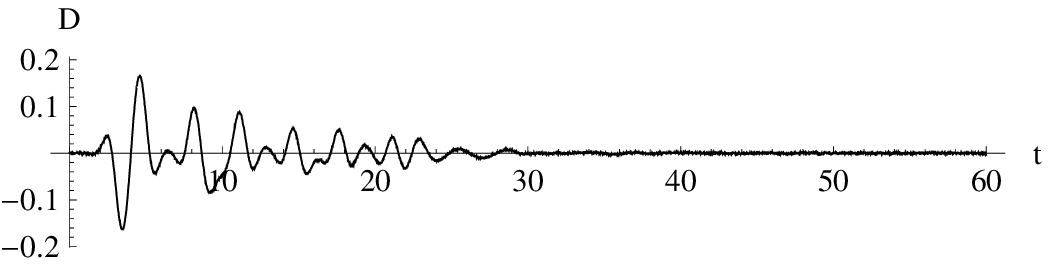}\\[10pt] 
\includegraphics[clip,trim=0in 0pt 0in 0in, width=1.5in]{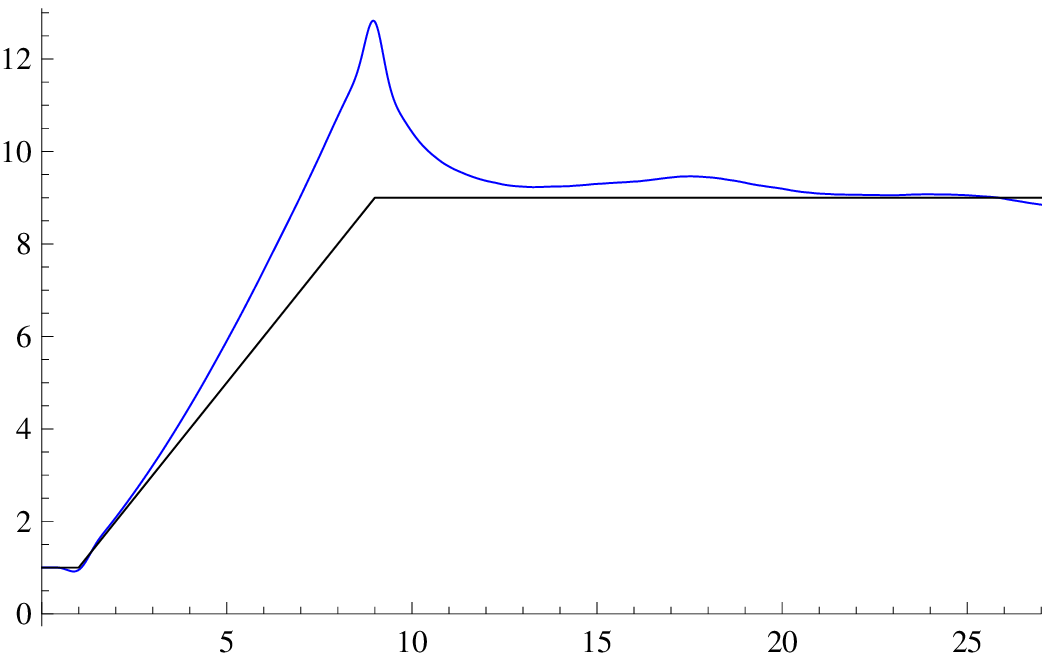}\hspace*{5pt}\includegraphics[clip,trim=0in 0pt 0in 0in, width=1.5in]{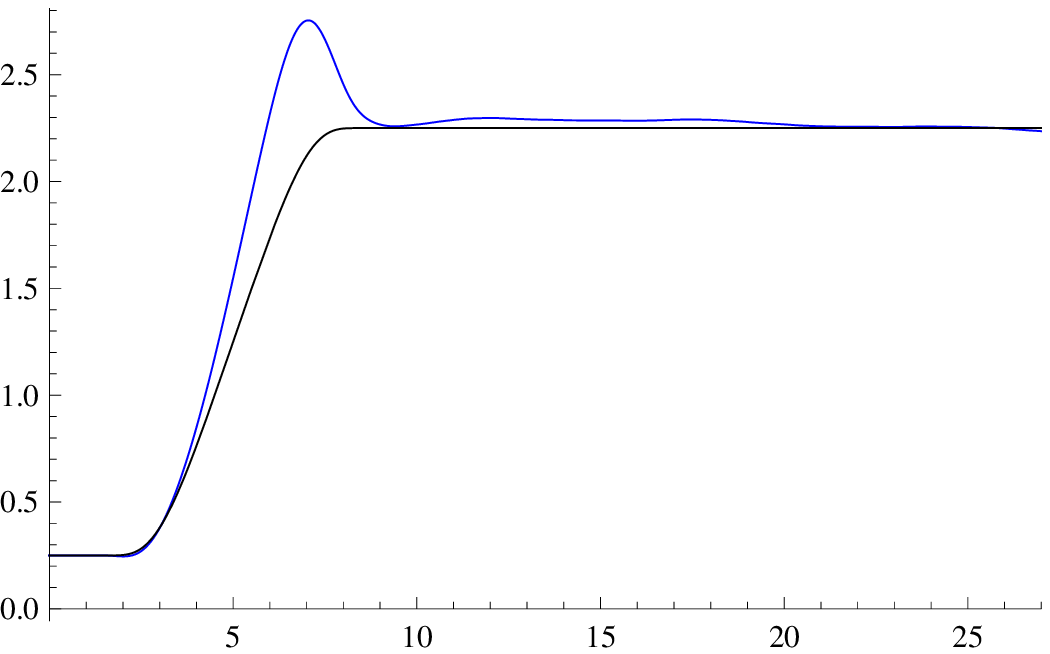}\hspace*{5pt}\includegraphics[clip,trim=0in 0pt 0in 0in, width=1.5in]{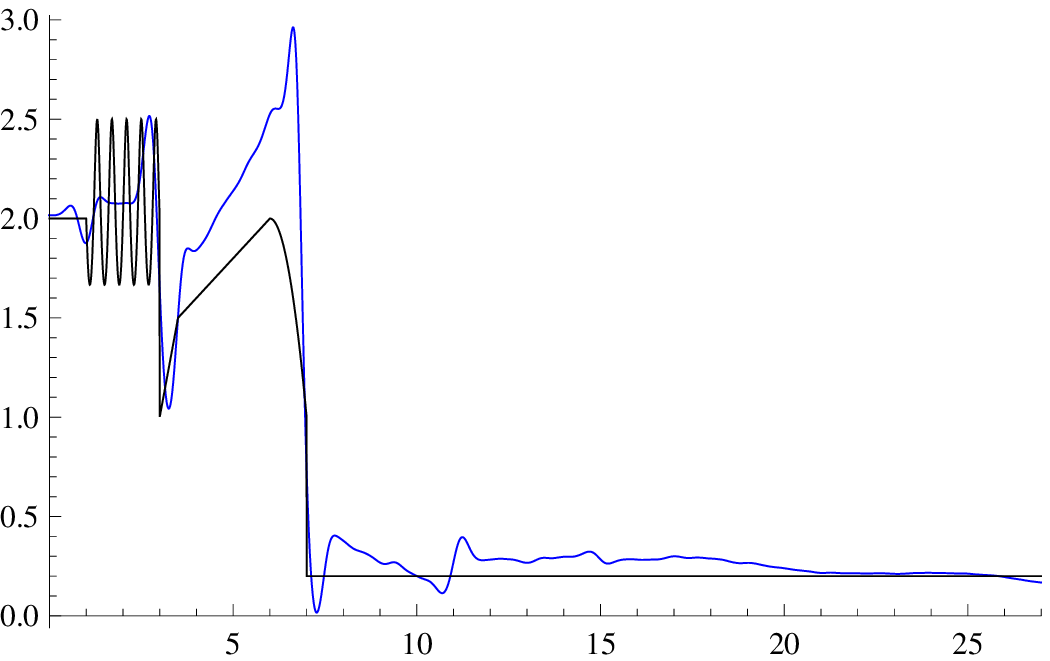}\hspace*{5pt}
\includegraphics[clip,trim=0in 0pt 0in 0in, width=1.5in]{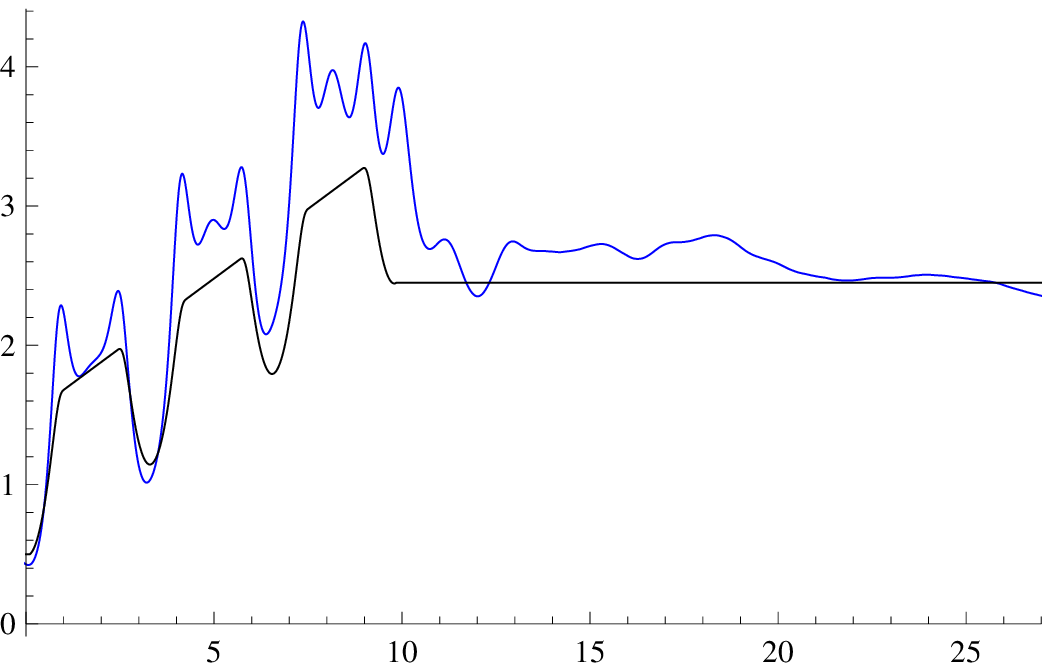}\caption{The effect of 5\% i.i.d.~Gaussian noise in the case of a zero mean source wave. Impedance profiles (blue, bottom row) reconstructed from noisy reflection data (top row) using the modified transform. The source wave is zero mean, plotted on the right of Figure~\ref{fig-wave-form}. }\label{fig-noisy-zero-mean}
}
}
\end{figure}

\section{Conclusion\label{sec-conclusion}}

The numerical results of the previous section illuminate two distinct applications of the refined impedance transform.  Firstly, there is the inverse problem of using reflection data to estimate the value $\zeta_+$ of the (eventually constant) impedance on the far side of an arbitrary, unknown slab supported on $[x_-,x_+]$.  This problem is solved directly by evaluating the refined impedance transform on recorded data $D$ or $D^{(-k)}$, as appropriate.  Moreover, the numerical examples show that the rate of convergence at play in Theorems~\ref{thm-slab} and \ref{thm-general} is reasonable---one does not require data much beyond the two way travel time to $x_+$.   By contrast, the standard approximation (\ref{standard}) does not converge to the correct value at all, producing relative errors that range from 25-94\%, even at large recording times.  
(Note that the problem of determining reflection data on the far side of unknown structure arises not only in seismic, but also in imaging biological tissue with coherent electromagnetic reflection data such as in \cite{WiFeWe:2006}.)

A second application of the refined impedance transform is to approximate $\zeta(x)$ for all $x<t_{\max}/2$.  The numerical examples illustrate that the refined impedance transform (or modified transform) gives a reasonable approximation to the full impedance profile, exhibiting its main qualitative features, albeit with overshoot and smoothing of details.  In certain cases, such as for the initial part of the third profile in Figure~\ref{fig-profiles}, the standard approximation (\ref{standard}) and the refined impedance transform give very similar results.  Roughly speaking, this is the case provided $\zeta(x)/\zeta_-$ is close to one.  If, on the other hand, $\zeta(x)/\zeta_-$ drifts outside the interval $(1/3,3)$, say, then substantial differences emerge between the standard estimate and the refined impedance transform---as happens in the latter part of the third profile in Figure~\ref{fig-profiles}---with the refined impedance transform tracking the value of the true impedance. Physical measurements of acoustic impedance of sedimentary rock taken from well bores show ratios for different rock types in excess of 4:1 (see \cite{Bo:1984}), and hence values of $\zeta(x)/\zeta_-$ substantially greater than 3. So in the context of seismic imaging the discrepancy between the refined impedance transform and the classical approximation of Peterson et al.~is important, with the former offering the prospect of significantly improved estimates.  

Iterative inversion methods in seismic (such as full-waveform inversion) require an initial estimate of acoustic impedance, the accuracy of which is crucial to convergence.  But there is no universally accepted way to obtain an initial estimate.  The refined impedance transform provides an objective approximation to impedance directly in terms of measured data---and thus may serve a useful role in the context of iterative methods.  

The refined impedance transform is relatively insensitive to noise because it involves the accumulation function of the data.  This is illustrated by comparison between Figure~\ref{fig-noisy-gaussian} with the earlier Figure~\ref{fig-wave}.  Ten percent noise does not seriously degrade the accuracy of the transformed data.  The situation is a little more delicate in the case where the source wave form has zero mean, however.  If $D$ includes noise, then the modified data $D^{(-k)}$ is perturbed by the $k$-fold integral of the noise, the variance of which increases with $k$.  Nevertheless, Figure~\ref{fig-noisy-zero-mean} illustrates that with $k=1$ and a noise level of 5\%, the refined impedance transform still gives reasonable results.  

In summary, the refined impedance transform comprises a new tool with which to analyze plane wave reflection data. It substantially outperforms the standard approximation while retaining the latter's speed and simplicity. 

%\newpage

\appendix

\section{Formulation in terms of pressure\label{sec-pressure}}

Depending on the context, it may be preferable to work with pressure $p(x,t)$ instead of particle velocity $u(x,t)$.  The present section compiles the formulas applicable to pressure that are analogous to those derived above for particle velocity.  Pressure is governed by the wave equation
\begin{subequations}\label{pressure-wave}
\begin{gather}
p_{tt}-\textstyle\zeta (\frac{1}{\zeta}p_x)_x=0\label{pressure-wave-equation}\\
p(x,0)=W(x)\qquad p_t(x,0)=-W^\prime(x)\label{pressure-initial-condition}.
\end{gather}
\end{subequations}
To avoid confusion with the previous notation let $K_\zeta$ denote the reflection Green's function for (\ref{pressure-wave}), and denote the associated measured data by 
\begin{equation}\label{pressure-data}
F=\widetilde{W}\ast K_\zeta.
\end{equation}
The classical estimate for $\zeta$ in terms of $K_\zeta$ is 
\begin{equation}\label{pressure-standard}
\zeta(x)\cong\zeta_-e^{2\int_{-\infty}^{2x}K_\zeta(t)\,dt}.
\end{equation}
The version of the refined impedance transform applicable to pressure data is the transform $\pimp_{w,c}$ defined by the formula
\begin{equation}\label{pressue-impedance-transform}
\pimp_{w,c}\;g(x)=\imp_{w,c}\bigl(-g\bigr)(x)=c\frac{w+\int_{-\infty}^{2x}g(t)\,dt}{w-\int_{-\infty}^{2x}g(t)\,dt}.
\end{equation}
It then follows from Theorem~\ref{thm-slab} that 
\begin{equation}\label{pressure-convergence}
\pimp_{w,\zeta_-}F(x)\rightarrow\zeta_+\qquad\mbox{ as }\qquad x\rightarrow\infty.
\end{equation}
The corresponding estimate for $\zeta$ in terms of the pressure Green's function is 
\begin{equation}\label{pressure-alternative}
\zeta(x)\cong\zeta_-\frac{1+\int_{-\infty}^{2x}K_\zeta}{1-\int_{-\infty}^{2x}K_{\zeta}}.  
\end{equation}

%\bibliography{References1.2017}

\end{document}